\documentclass[12pt]{article}
\usepackage{microtype} %% Microtypographical adjustments.  If you do                       %% not know what it does, do not worry. :-)

\usepackage{graphicx}
\usepackage{amsmath}%
\usepackage{mathtools} % includes amsmath and a few more gadgets
\usepackage{amsthm}
\usepackage{amsfonts}
\usepackage{amssymb}
\usepackage[margin=0.75in]{geometry}
\usepackage{color}
\usepackage{latexsym}
\usepackage{algorithm}
\usepackage[noend]{algorithmic}
\usepackage{subfigure}
\usepackage{soul}
\usepackage{enumerate}
\usepackage{multirow}
\usepackage[nameinlink]{cleveref}
\usepackage{thm-restate}
\usepackage{xfrac}

\usepackage{diagbox}
\usepackage{makecell}
\usepackage{array,booktabs}
\usepackage{hhline}
\usepackage[table]{xcolor}
\newcommand\VRule[1][\arrayrulewidth]{\vrule width #1}

\usepackage{comment}%for multi line comments

\newtheorem{theorem}{Theorem}%[section]

\newtheorem{theoremB}{Theorem}%See comment on next line.
\newtheorem{lemma}[theorem]{Lemma}

\theoremstyle{remark}

\theoremstyle{definition}

\newtheorem{remark}[theorem]{Remark}

\newenvironment{Proof}[1]{\par\noindent{\bf Proof{#1}:}\hspace{3pt}}{}
\newcommand{\qedForProof}{\qed\vspace{10pt}} %% add a  square  and 10 pt vspace at the end of the Proof. 

\newcommand{\old}[1]{{{}}}

\graphicspath{{./Figures/}}

\def\N{{N}}
\def\R{\textbf{R}}

\let\ll\lambda

\def\thA{\theta_1}
\def\thB{\theta_2}
\def\aA{2}
\def\aB{2}
\def\el{e^{\lambda}}
\def\ec{e^{-c}}

\def\ql{\left(1-{\ll}/{n}\right)}
\def\an{{n-an/\log n}}
\def\anlln{{n-an\log \log n/\log n}}

\def\anllnq{{an(\log \log n-1)/\log n}}
\def\la{\lambda a}

\def\lnln{\log\log n}
\def\lnlnln{\log\log\log n}
\def\lln2{\log^2\log n}
\def\integrandbase{\left({\ec\log n}/{n}\right)}
\def\integrallimitL{\an}
\def\Exp{\textnormal{Exp}}
\def\RHS{right-hand side }
\def\LHS{left-hand side }

\def\errLognly{O\left({1}/{\log n}\right)}
\def\errLognst{O\left(\tfrac{1}{\log n}\right)}
\def\errOne{O\left(1\right)}

\def\lglgFlgly{\log\log n/\log n}
\def\lglgFlgst{\frac{\log\log n}{\log n}}
\def\lglgFlgtst{\tfrac{\log\log n}{\log n}}
\def\lemmafunctionG{n\lnln/(\ll\log n)}
\newcommand{\sizedSubNumber}[2]{{#1}_{\scalebox{0.7}{{#2}}}}
\newcommand{\sizedSubLetter}[2]{{#1}_{\scalebox{0.7}{\textit{#2}}}}
\newcommand\Item[1][]{%
  \ifx\relax#1\relax  \item \else \item[#1] \fi
  \abovedisplayskip=0pt\abovedisplayshortskip=0pt~\vspace*{-\baselineskip}}%Dina: to align "item" with equations. If no label, try option t of aligned: $\! \begin{aligned}[t] equation.... end...$ 

\usepackage{xspace}
\usepackage[numbers,sort&compress]{natbib}%sorts citations by number and also compresses [1,2,3] to [1--3]

\newcommand{\Eq}[1]{(\ref{#1})}

\def\Thanks#1{\gdef\thefootnote{\arabic{footnote}}\thanks{#1}}
\def\ThanksComma#1{\gdef\thefootnote{\arabic{footnote},}\thanks{#1}{
}}

\begin{document}
\title{The Time for Reconstructing the Attack Graph in DDoS Attacks}

\author{D. Barak-Pelleg
%\ThanksComma{Corresponding author. Tel.: +972 50 4271414; E-mail addresses: dinabar@post.bgu.ac.il; Current  address: dina.barak.pelleg@gmail.com (D. Barak-Pelleg)}
\ThanksComma{Department of Mathematics, Ben-Gurion University, Beer Sheva 84105, Israel. E-mail: dinabar@post.bgu.ac.il; Current  address: dina.barak.pelleg@gmail.com}
\Thanks{Research supported in part by Cyber Security Research Center, Prime Minister's Office, and a Hillel Gauchman scholarship.} 
\and  
D.~Berend
\ThanksComma{Departments of Mathematics and Computer Science, Ben-Gurion University, Beer Sheva 84105, Israel. E-mail: berend@math.bgu.ac.il}
\ThanksComma{Research supported in part by the Milken Families Foundation Chair in Mathematics and Cyber Security Research Center, Prime Minister's Office.}
\Thanks{Research supported in part by the Center for Advanced Studies in Mathematics at Ben-Gurion University.}
}

\date{}
%\tableofcontents\newpage

\maketitle

\begin{abstract}
Despite their frequency, denial-of-service (DoS) and distributed-denial-of-service (DDoS) attacks are difficult to prevent and trace, thus posing a constant threat.
One of the main defense techniques is to identify the source of attack by reconstructing the attack graph, and then filter the messages arriving from this source. One of the most common methods for reconstructing the attack graph is Probabilistic Packet Marking (PPM). We focus on edge-sampling, which is the most common method.  
Here, we study the time, in terms of the  number of packets, the victim needs to reconstruct the attack graph when there is a single attacker. This random variable plays an important role in the reconstruction algorithm.	Our main result is a determination of the asymptotic distribution and expected value of this time. 
	
The process of reconstructing the attack graph is analogous to a version of the well-known coupon collector's problem (with coupons having distinct probabilities). Thus, the results may be used in other applications of this problem.


\end{abstract}

\section{Introduction}\label{introduction}

\subsection{DDoS Attack and PPM}
A \textit{denial-of-service} (DoS) attack is a cyber attack in which the victim, a particular computer on the internet network, is assailed by a single attacker, seeking to make the victim unavailable for service. This goal is accomplished by flooding the victim with fake data packets until it is unable to fulfill legitimate requests, or even collapses. A distributed-denial-of-service (DDoS) attack is similar, but with multiple attackers. Both types of attack are common as they are quite easy to launch. Despite their frequency, these attacks are difficult to prevent and trace, thus posing a constant threat (see \citep{DDoSNews} for the latest DDoS attack news). 

Several defense techniques and tools are available to deal with these attacks; usually, a combination of  approaches is employed  (see, for example, \citep{LoukasOke,ZargarJoshiTipper} for surveys on defense techniques).  One of the main approaches is to identify the source of attack, and then filter the messages arriving from this source. There are a few methods to implement  this approach~\citep{BelenkyAnsari}.  One of these methods is by reconstructing the attack graph. This graph is a tree type graph, in which the root represents the victim, the leaves represent the attackers, and the internal nodes  represent the routers connecting the attackers to the victim. (Thus, in a DoS attack, the graph comprises a path.) There are various methods for reconstructing the attack graph (see \citep{Kiremire2014}). 
In the current work we focus on Probabilistic Packet Marking (PPM), introduced in \citep{BurchCheswick}. Specifically, we deal with edge-sampling, the most common method used in PPM.

In edge-sampling, there are two processes taking place simultaneously. The first is on the routers side: Each router in the network, upon receiving a packet, and before forwarding it, decides at random whether to mark it or not; The marking probability is~$p$ (fixed for all routers). If the packet has already been marked by a previous router, the new mark will override the old one. 
Thus, the probability of a packet received by the victim to carry the mark of the router at distance $i$ from him is \(p(1-p)^{i-1}\). When a router marks a packet, it writes there its identity, and the next router (if it does not override the mark) adds to it its own identity and starts a counter. When any router farther along the path decides not to override the mark, it increases the counter by~$1$. Thus, when the victim receives a marked packet, the mark consists of the edge in the attack path corresponding to the (last) marking router and the router following it, and the distance of this edge from the victim.
The second process is on the victim's side: The victim collects the marks in order to reconstruct the attack graph.  

The victim starts collecting  marks upon suspecting he is under attack; that is, when there is a sudden jump in the  arrival rate of packets. When should this process be terminated? Namely, when should the victim  decide it has obtained enough data in order to reconstruct the full attack graph? On the one hand, the longer the victim continues collecting marks, the greater the chance of being able to reconstruct the full attack graph. On the other hand, if the victim waits too long, it might collapse by the flood of incoming packets.
The time (in terms of the  number of packets) the victim needs  in order to reconstruct the full attack graph, when there is  a single attacker, is also referred to as the Completion Condition Number  \citep{Sairam-Saurabh-2013}. This random variable, which we will denote by $D$, plays an  important role in the reconstruction algorithm. 
Savage, Wetherall, Karlin, and Anderson \citep{AKSW} considered the  expected number $E(D)$ of packets needed, and showed that
$E(D)\le {\ln{n}}/({p(1-p)^{n-1}})$, where $n$ is the distance of the attacker from the victim. Thus, they suggested to wait until obtaining ${\ln{n}}/({p(1-p)^{n-1}})$  packets. 
Sairam and Saurabh \citep{Sairam-Saurabh-2013} showed that, in many cases, this number of packets may not be enough. They up-bounded the standard deviation of $D$ and suggested to add a third of this bound to the above bound on $E(D)$, thus increasing the reliability of the algorithm.

The process of obtaining the marks by  the victim is analogous to a version of the coupon collector's problem \citep{AKSW,Sairam-Saurabh-2013,Sairam-Saurabh-2016,Shioda}. We now recall this classical problem.

\subsection{The Coupon Collector's Problem}
Suppose that a company distributes packages of some product and that each package contains a single coupon. There are $n$ types of coupons, and a customer wants to collect them all. Each time that he buys a package, he gets one of the types uniformly at random. We want to know how many packages need to be purchased on the average until getting all types of coupons.  The problem goes back at least as far as de Moivre,  who mentioned it in a collection of problems regarding various games of chance \cite{Moiver}.
The solution to this problem has been known for many years; the expected number of coupons we need to draw is $nH_n$,
where $H_n=1+{1}/{2}+{1}/{3}+\cdots +{1}/{n}$ is the $n$-th harmonic number.
Asymptotically, this expectation is $n\ln n+\gamma n+O(1)$, where $\gamma = 0.577\ldots$  is the Euler-Mascheroni constant.

The problem, and various extensions thereof, have drawn much attention for many years (see, for example,  \citep{BarakEtAl2022A,Flatto,MyersWilf2006,Laplace,ErdosRenyi1961,NewmanShepp,KobzaJacobsonVaughan}; see also the surveys \cite{CouponHistory,Dawkins}). One of the extensions, considered by von~Schelling \citep{VonSchelling1954}, and Flajolet,  Gardy,  and  Thimonier \cite{FlajoletCoupon}, dealing with the case where various coupons show up with distinct probabilities, turns out to be very relevant to our problem. 
In the next subsection we will see that the reconstruction of the attack graph is naturally translated to this variant. 

\subsection{Edge-Sampling and  Coupon Collecting}\label{introduction1}
As mentioned above, in a DoS attack, the attack graph is just a path. Denote its  
vertices by $v_0,\ldots,v_n$, where  $v_0$ represents the victim and $v_n$ represents the attacker, and  its edges by $e_i=\{v_{i-1}, v_i\}$, $1\le i\le n$. Each $e_i$  represents the link between the router at distance $i-1$ with that at distance $i$ from the victim. 

To connect the reconstruction problem with the coupon collector's problem, we regard the victim of the DoS attack as a coupon collector, and each $e_i$ as the $i$-th type coupon. The event ``the victim has obtained a packet marked by the link at distance of $i-1$ from him" is translated to ``the coupon collector has received a coupon of type $i$".  Obtaining the marks of all links of the attack path is equivalent to the collector having obtained all coupon types.

As indicated above, the version of the coupon collector's problem we have here is where the coupons have distinct probabilities. 
Each coupon type $i$ is drawn with probability $p_i=p(1-p)^{i-1}$. Note that the sum of these probabilities is $\sum_{i=1}^n p_i = 1-(1-p)^n<1$, as at each step there is a probability of $(1-p)^n$ to obtain an unmarked packet. Thus, it will be convenient for us to add a ``dummy" coupon of type $0$, whose probability is $p_0=(1-p)^n$, and a corresponding ``dummy" edge $e_0$ to the attack path. This addition is inconsequential for the following reason. We take the marking probability to be $p=\lambda/n$, for some arbitrary fixed $\lambda>0$, and assume that $n$ is large. Hence all ``real" coupons have probabilities $\Theta(1/n)$, while the probability of the dummy coupon is $\Theta(1)$. The probability for the dummy coupon to be obtained last is therefore extremely small. Whether the goal is to collect only all real coupons, or it is to collect also the dummy one, is immaterial; the dummy coupon will anyway (most probably) arrive long before all real coupons have arrived.

\subsection{Paper Organization}
In Section \ref{sec:Results} we define a continuous analogue of our problem, which is more convenient to deal with than our discrete model. Next we state the main results, first for the continuous version, and then for the discrete one. We note that the convergence rate in our theorems is quite slow. Thus, Section \ref{simulations} describes simulations we preformed for both models; the simulations hint that, indeed, the convergence rate is not much faster than what is guaranteed by the main results. In Section \ref{Proofs} we prove the results for the continuous model, and then explain how they can be used to prove those on the original model. 
%In Section \ref{ProofByNeal} we show that, trying to use the results of \cite{Neal2008}, we would be able to prove only some of our results, and with no simplification of the proofs.

\bigskip

\bigskip

\section{Main Results}\label{sec:Results}
%aaaaaaaaaaaaaaaaaaaaaaaaaaaaaaa

Let us first consider a continuous version of our problem. 
    The idea of using a continuous model has been used several times in the classical case  (see  \citep{Holst1986,Ilienko2020}). In this model  there are~$n$ independent, incoming flows of coupons
$$\sizedSubNumber{T}{1}\sim\Exp\left(\sizedSubNumber{p}{1}\right),\ldots,\sizedSubLetter{T}{n}\sim\Exp\left(\sizedSubLetter{p}{n}\right),$$
where $\sizedSubLetter{T}{i}$ is the inter-arrival time between  consecutive coupons of type $i$. Same as in the regular model, we are interested in the waiting time until all  coupon types arrive. 
Differently from the regular model, the waiting times are exponential  instead of geometric. Also, in the continuous model the variables are independent, whereas in the discrete model they are not. 
Thus, the probability that  the $i$-th coupon type  has not been seen until time $t$ is
$$e^{-p_it}=e^{-\ll/n\ql^{i-1}t}.
$$
%In particular, the probability we have not seen the first coupon until time $t$ and the last are:$$e^{-\ll/n\ql^{n-1}t},\qquad e^{-\ll/n \cdot t},$$respectively.
%By independence, the probability we have not seen coupon $i$ and also coupon $j$ until time $t$ is$$e^{-p_it}\cdot e^{-p_jt}=e^{-(p_i+p_j)t}.$$
Denote by $T$ the time until we get all coupons:
$$T=\max_{1\le i\le n}\sizedSubLetter{T}{i}.
$$
Given a sequence $\left(X_{n}\right)_{n=1}^{\infty}$ of random variables
and a probability law $\mathcal{L},$  write~$X_{n}\xrightarrow[n\to\infty]{\mathcal{D}}\mathcal{L}$
if the sequence converges to $\mathcal{L}$ in distribution.
Recall  that   a random variable $X$ is \textit{Gumbel distributed} with parameters $\mu\in \bf R$ and $\beta>0$,  and we write $X\sim\textnormal{Gumbel}(\mu,\beta)$, if its distribution function is given by  \citep{Gumbel1954,PinheiroFerrari}: %in the book it is on page 21: 3.6 and 3.7.
\begin{equation}\label{gumbel}
    F(x)=e^{-e^{{-(x-\mu)}/{\beta}}},\qquad x\in \bf R.
\end{equation}
\begin{theorem}\label{ConvInDist}
%\begin{description}
%\item[{ $a.$}]
%location in main4:Thm \ref{ConvInDist}
The asymptotic distribution of the waiting time for all coupons in the continuous model is given by:
\begin{equation*}
    \frac{T-(\el/\ll)\cdot n(\log n-\log \log n)}{n}
    \xrightarrow[n\to\infty]{\mathcal{D}}\textnormal{Gumbel}\left(-\tfrac{\el}{\ll}\log\ll,\tfrac{\el}{\ll}\right).
\end{equation*}
%\item[{ $b.$}]
%\end{description}
\end{theorem}
We will actually prove the following stronger version of the theorem, which provides information about the rate of convergence in Theorem \ref{ConvInDist}.
Denote:
\begin{equation}\label{T'notation}
    T'=\frac{T-(\el/\ll)\cdot n(\log n-\log \log n%-\log\ll
    )}{%(\el/\ll)
    n},\qquad
\end{equation}
\setcounter{theoremB}{0}
\begin{theoremB}\label{ConvInDist'}
%location in main4: Thm \ref{F_T'}
For $t'\in\bf{R}$ and $n\to\infty$,
\begin{equation*}
    {F}_{T'}(t')=\exp\left(-e^{-\left(t'-\left(-\el \log\ll/\ll\right)\right)/\left(\el/\ll\right)}\right)+O\left({\log\log n}/{\log n}\right).
\end{equation*}
\end{theoremB}
Getting back to the discrete model, recall that $D$ is the number of coupons we need to  collect in order to get all real types in the discrete case.
Similarly to (\ref{T'notation}), denote:
\begin{equation*}
    D'=\frac{D-(\el/\ll)\cdot n(\log n-\log \log n%-\log\ll
    )}{%(\el/\ll)
    n}.
\end{equation*}
\begin{theorem}\label{F_D}
%\begin{description}
%\item[ $a.$]%location in main4: first part of Thm \ref{F_D'}.a, or equation number
The asymptotic distribution  of the time required for reconstructing the attack graph are given by:
\begin{equation}\label{F_D'.1}
    D'%\frac{D-(\el/\ll)\cdot n(\log n-\log \log n%-\log\ll
    %)}{%(\el/\ll)
    %n} 
    \xrightarrow[n\to\infty]{\mathcal{D}}
         \textnormal{Gumbel}\left(-\tfrac{\el}{\ll}\log\ll,\tfrac{\el}{\ll}\right).
\end{equation}
%\item[ $b.$]%location in main4: second part of Thm \ref{F_D'}.a, or equation number
Moreover, as $n\to\infty$,
\begin{equation}\label{F_D'.2}
    {F}_{D'}(d')=\exp\left(-e^{-\left(d'-\left(-\el \log\ll/\ll\right)\right)/\left(\el/\ll\right)}\right)+O\left({\log\log n}/{\log n}\right),\qquad d'\in\bf{R}.
\end{equation}
%\end{description}
\end{theorem}

Note that the convergence rate we obtain is rather slow, which goes hand in hand with the rate of convergence of other quantities related to the coupon collector's problem \cite{Brayton63,Holst1986}.
In the next section we describe a large simulation we have performed, which hints that the error term is probably near-optimal.

\begin{theorem}\label{Expectation}
\begin{description}
\item[ $a.$]%location in main4: Thm \ref{F_D'}.b
The expected times until we get all coupons in the two models coincide:
\begin{equation*}
    E(D)= E(T).%=\tfrac{e^{\lambda}}{\lambda}\cdot n\left(\log n-\log \log n+\gamma-\log\lambda\right)+O\left({n\log\log n}/{\log n}\right).
\end{equation*}
\item[ $b.$]%location in main4: Thm \ref{Expectation}
As $n\to\infty$:
\begin{equation}\label{eq-expectation}
    \begin{aligned}
    E(T)= \tfrac{e^{\lambda}}{\lambda}\cdot n\left(\log n-\log \log n+\gamma-\log\lambda\right)+O\left({n\log\log n}/{\log n}\right).
\end{aligned}
\end{equation}
%where $\gamma$ is the Euler-Mascheroni constant.
\end{description}
\end{theorem}
\begin{remark}\label{remarkNeal}
   In principle, we could have used the results in \citep{Neal2008} to prove Theorem~\ref{F_D}. However, this would lead to the same type of calculations. More importantly, the estimates we would have received would not be strong enough to prove Theorem~\ref{Expectation}.
\end{remark}

\begin{figure}[ht]
\centering
\includegraphics[width=0.5\textwidth]{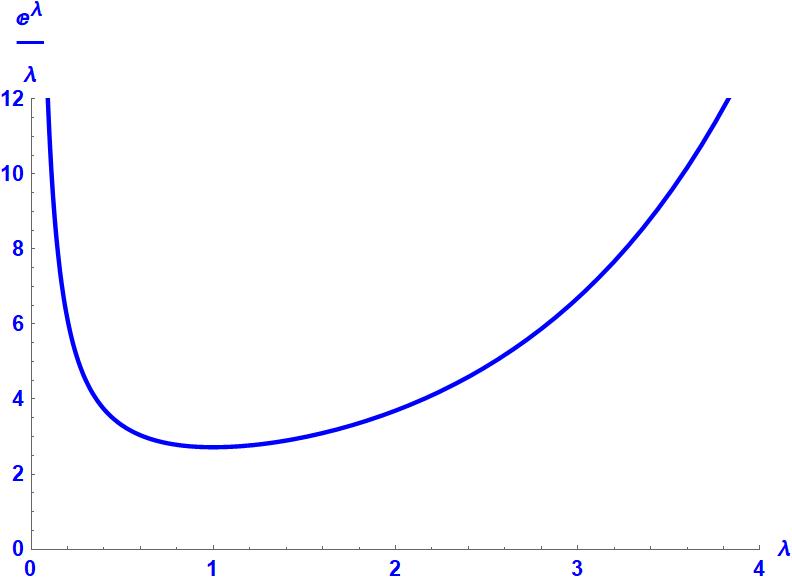}
\caption{The effect of the coefficient $\ll$ on the reconstruction time}% $e^\ll/\ll$}
\label{fig1}
\end{figure}

\begin{remark}
According to the theorem, the reconstruction time is roughly proportional to $e^\ll/\ll$. For $\ll>0$, the expression $e^\ll/\ll$ is minimal at $\ll =1$ (see Figure \ref{fig1}). Hence, as $n\to\infty$, the expectation $E(D)$ will be minimal very close to the point $\ll =1$. Thus, the optimal choice for the edge-sampling algorithm is $p=1/n$ (as claimed by Savage et al.\ \citep[p.300]{AKSW}). Thus, we have held our simulations only for $\lambda=1$.
\end{remark}

\bigskip

\bigskip

%\newpage

\section{Simulation Results}\label{simulations}
We have performed a simulation for the time needed to collect all types of (real) coupons.
As mentioned above, the convergence of quantities in CCP is rather slow. The results of the following  simulations  hint already    that the order of magnitude of the error obtained in  Theorem~\ref{Expectation} %**
is close to optimal. 

 In our experiments,  $\ll=1$, $n= 10^4$, and the number of iterations of each test is $M=10^5$. Everything has been performed on Mathematica, and we point out several technical points that may be on interest to its users.

The simulation was preformed for both the discrete model and the continuous one.
In the discrete case, in each of the $M$ runs we have drawn the coupons one by one, each drawing being independent of the others. In each drawing, the coupon of type $i$ was selected with  probability $p_i=1/n(1-1/n)^{i-1}, 1\le i\le n$. We have continued the process until all $n$ types of real coupons have been drawn and saved the number of drawings. Thus, we have obtained a list of length $M$, of the times at which the various iterations completed their runs.

In the continuous case, in each of the $M$ iterations we selected $n$ random exponential
variates with parameters $p_1,\ldots ,p_n$, and took their maximum.

In Table \ref{table1}, the first two columns present the sample means (rounded to the nearest integer) received in the two experiments. The third column shows the main term $e n\left(\log n-\log \log n+\gamma\right)$ on the right-hand side in our expression for $E(T)$ and $E(D)$ from Theorem~\ref{Expectation}.  %**
The last column presents the order of magnitude of the error term, namely $n\log\log n/\log n$. Note that the two means are relatively very close, and both are in line with the theoretical main term, given the allowed error. Thus, the error term in Theorem~\ref{Expectation} %**
may well be of the correct order of magnitude.
\begin{centering}
\setcellgapes{7pt}
\begin{table}[ht]
\makegapedcells
\centering
\begin{tabular}{!{\VRule}c!{\VRule}c!{\VRule}c!{\VRule}c!{\VRule}}
\hline%\specialrule{2pt}{0pt}{0pt}%thick line
%\multicolumn{1}{!{\VRule}c!{\VRule}}{$m$}

\multicolumn{1}{!{\VRule}c!{\VRule}} {$\overline{D}_M$} &
 $\overline{T}_M$ & $e n\left(\log n-\log \log n+\gamma\right)$ & ${n\log\log n}/{\log n}$ \\ \hline
& & & \\[-2.6ex] \hline
%First line11111111111111111111111
$207945$ & $207885$  & $205699$ & $2410$\\ \hline
 %Second line222222222222222222222
%$2$& $0.213$ & $0.167$  & $0.225$\\ \hline
% \hline%\specialrule{2pt}{0pt}{0pt}%thick line
\end{tabular}
\caption{The sample means vs.\ the theoretical results on the expectation.}
\label{table1}
\end{table} 
\end{centering}

\begin{figure}[ht]
\centering
\includegraphics[width=0.75\textwidth]{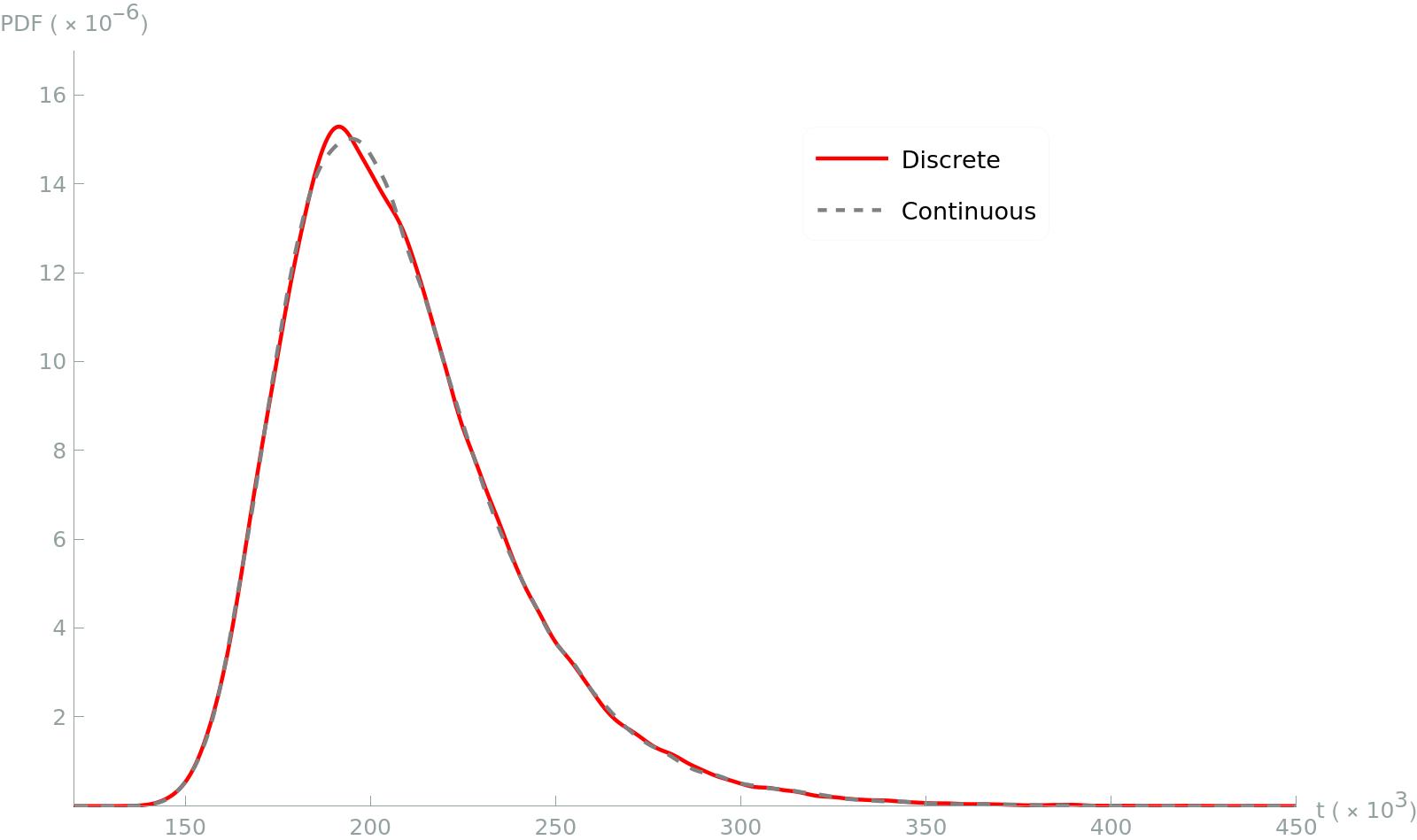}
\caption{Smoothed histogram for the results of both models.}
\label{plot1}
\end{figure}

Not only the sample means are close, as may be seen in Table \ref{table1}. In Figure \ref{plot1} we present the (smoothed) PDFs of the simulation data for both models (using the default option  \textquotesingle\textquotesingle PDF\textquotesingle\textquotesingle\  in Mathematica's SmoothHistogram). The results for the discrete model  presented by the smooth red line and those of the continuous by the dashed gray line.

We have utilised the function $\textnormal{FindDistribution}$ to find   Mathematica's guess for the most fitting distribution for the sample of $D$. We have repeated the simulation several times. In most cases, Mathematica guessed that the sample data is from a $\textnormal{Gumbel}$ distribution, but this was not always the case. In the simulation we have presented here, we  received three guesses: 

(i) First, we have given $\textnormal{FindDistribution}$ only the data of the simulation, without any ``hints" as to the required distribution. In this case, Mathematica guessed that the sample from $D$ is   from a mixture  of two basic distributions: $$0.73\times\textnormal{Gamma}(101, 1934)+0.27\times\textnormal{LogNormal}(12.33, 0.19).$$ 

(ii) Second, we have added the option $\textnormal{MaxItems}\to 1$, which yields a single,  best fitting distribution for the data. In this case, Mathematica's guess was  $\textnormal{Gubmel}(192983, 25762)$. 

(iii) We have noticed that the mean and variance of the distribution suggested in the second guess do not fit those of our sample. Thus, we have specified for Mathematica to find the most fitting $\textnormal{Gubmel}$ distribution by utilising the option $\textnormal{TargetFunctions}\to\lbrace{ \textnormal{ExtremeValueDistribution}\rbrace}$ (which is the name Mathematica uses for the distribution called Gumbel in our paper). In this case, Mathematica's guess was $\textnormal{Gubmel}(193800, 24506)$. 

\begin{figure}[ht]
\centering
\includegraphics[width=0.75\textwidth]{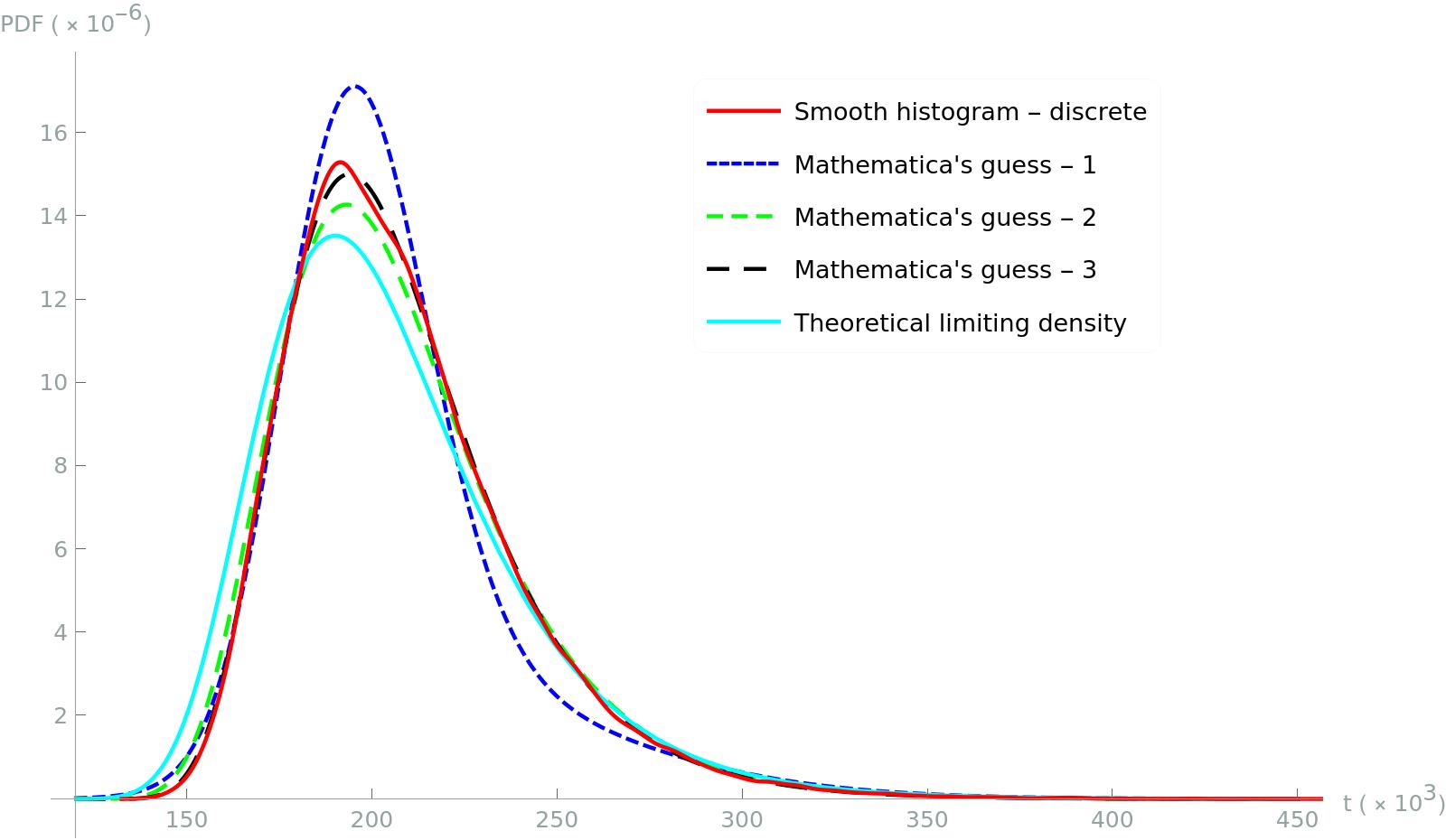}
\caption{Smoothed histogram of simulation vs. Mathematica's guesses and the theoretical limiting distribution PDFs.}
\label{plot2}
\end{figure}

 In Figure \ref{plot2} we present five graphs, generated by Mathematica. Four of them are based on the simulation data for the discrete model and the last depicts the prediction of the theoretical result. The red continuous line presents the (smoothed) probability density function of the simulation data, same as in Figure \ref{plot1}. The three dashed lines present the probability density function  of the three guesses (i)-(iii) above of Mathematica for the distribution most fitting the simulation data. The first is presented by a blue line of small dashes, the second -- by a green line of medium dashes, and the third -- by a black line of large dashes.
 The solid cyan line presents the density function of the $\textnormal{Gumbel}(en\log n-en\log\log n,en)=\textnormal{Gumbel}(190008,27183)$ distribution. This distribution is the approximation of the distribution of $D$, corresponding to the approximation of the distribution of $D'$ by $\textnormal{Gumbel}(0,e)$, as in (\ref{F_D'.1}).

%The solid cyan line presents the density function of the asymptotic distribution on the right-hand side of (\ref{F_D'.2}) (without the error term). This result  is given for $D'$, and therefore is scaled differently than the first four graphs. In order  to fix this scaling we need to go back from $D'$ to $D$. More specifically, we take the density  of the distribution function -- $F_{}\left((t - en\log n + en \log\log n)/n\right)$, where $F(t)=\exp(-e^{-t/e})$ is the distribution function of  $\textnormal{Gumbel}(0,e)$. Namely, the fifth graph is of $\tfrac1n f_{}\left((t - en\log n + en \log\log n)/n\right)$, where $f(t)=\tfrac{1}{e}e^{-t/e}\cdot\exp(-e^{-t/e})$.

\begin{figure}[ht]
\centering
\includegraphics[width=0.75\textwidth]{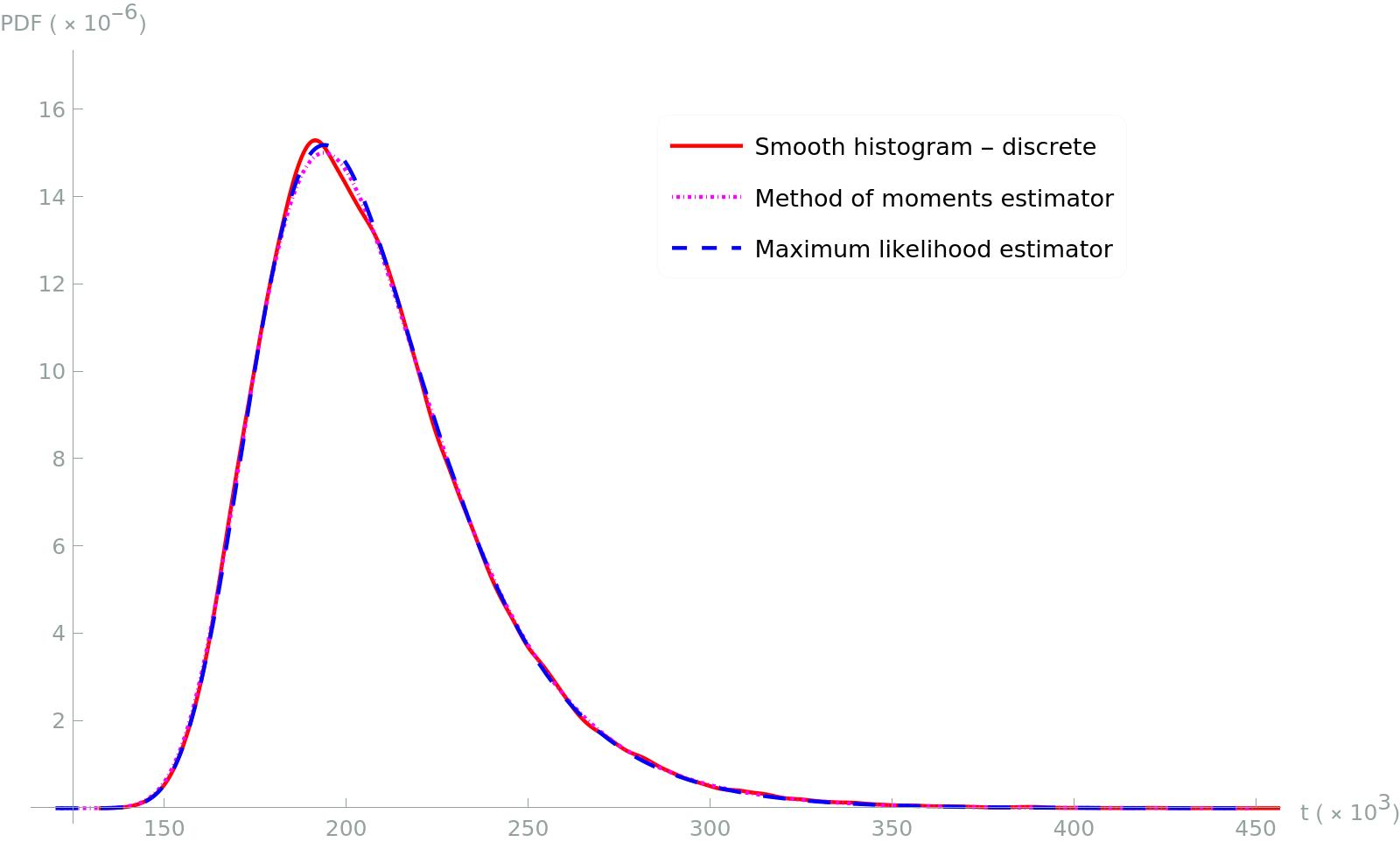}
\caption{Smoothed histogram of simulation data vs.\ Gumbel PDFs with parameters estimated by the method of moments and maximum likelihood.}
\label{plot2a}
\end{figure}

Note that here, when providing Mathematica with the hypothesized distribution type, Gumbel, we have an estimation problem of two unknown parameters $\mu$ and $\beta$. The simplest way to estimate these parameters is by the method of moments \citep{MahdiCenac2005,Corsini1995}. Employing Mathematica's $\textnormal{EstimatedDistribution}$ with the option $\textnormal{ParameterEstimator}\to$ \textquotesingle\textquotesingle MethodOfMoments\textquotesingle\textquotesingle, we get the same parameters as  guess (iii) above. Recall that the  method of moments estimator employs the sample moments to estimate the parameters.  Thus, as expected, in this case we get a Gumbel distribution whose expectation and variance fit the sample mean and sample variance. For maximum likelihood estimation  \citep{Fiorentino1984,Kimball1946}, the parameters are given implicitly, and thus more difficult to obtain. Employing Mathematica's $\textnormal{EstimatedDistribution}$ with the option $\textnormal{ParameterEstimator}\to$ \textquotesingle\textquotesingle\textnormal MaximumLikelihood\textquotesingle\textquotesingle,  we get $\textnormal{Gubmel}(193878, 24218)$. 
 In Figure \ref{plot2a} we depict three graphs, generated by Mathematica. As in Figure \ref{plot2}, the red continuous line represents the (smoothed) probability density function of the simulation data. The magenta dotted line is of the Gumbel distribution whose parameters were estimated by the method of moments, and the blue dashed line -- for the maximum likelihood estimator.
 
\begin{figure}[ht]
\centering
\includegraphics[width=0.75\textwidth]{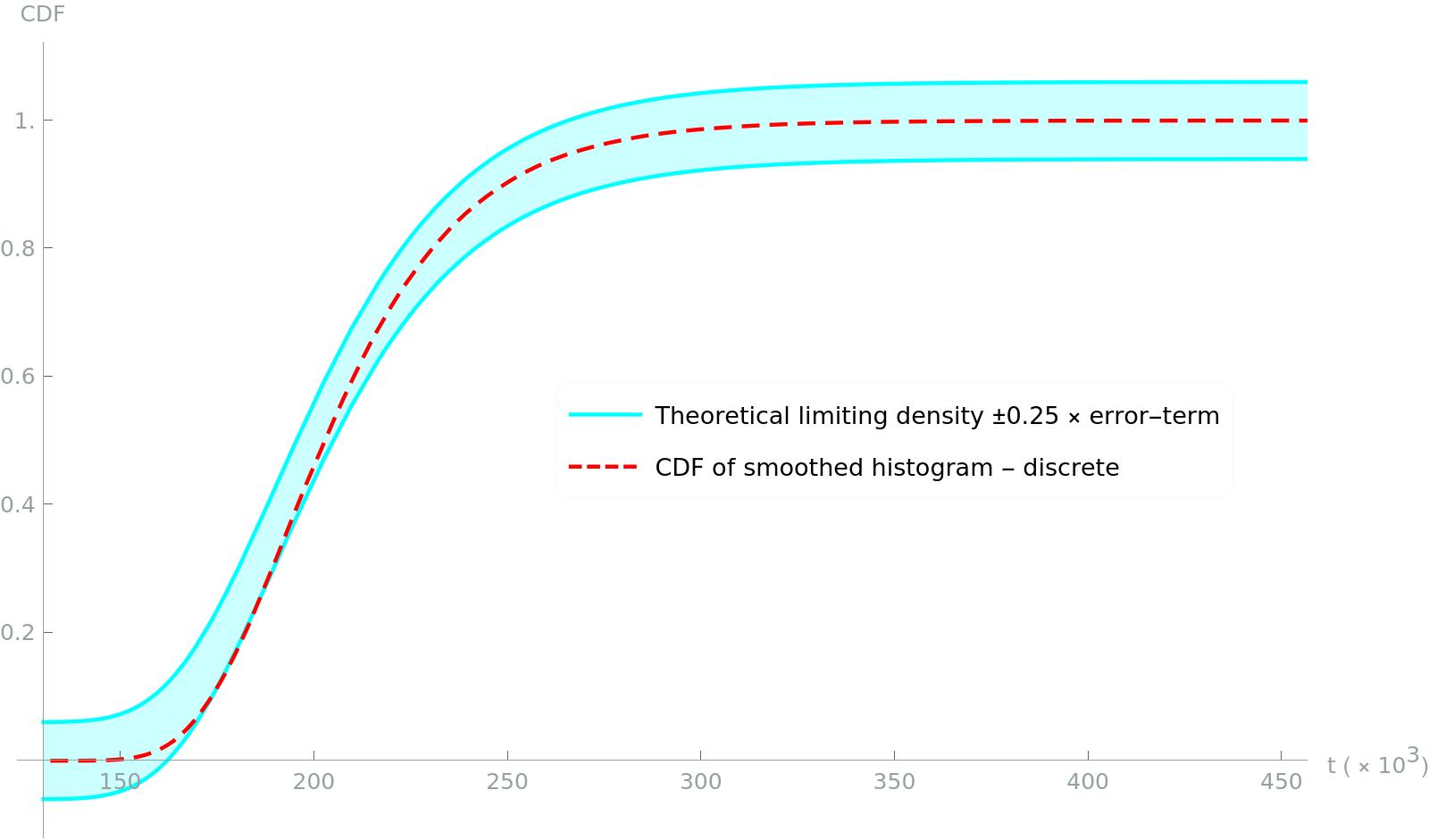}
\caption{The CDF of the theoretical limiting distribution +/- twice the error term vs. the CDF of the smoothed histogram of the simulation.}
\label{plot3}
\end{figure}

In Figure \ref{plot3} we have three graphs. The cyan solid lines both represent the main term on the \RHS of (\ref{F_D'.2}), raised and lowered by $0.25\log\log n/\log n$, namely one fourth of the expression in the error term. Explicitly, the graphs are of the functions 
$$F\left((t - en\log n + en \log\log n)/n\right) \pm 0.25\log\log n/\log n,$$
where $F(t)=\exp(-e^{-t/e}).$ 
We illustrate the closeness of the simulation results to the theoretical result by adding the CDF of the smoothed histogram of the data (using the option \textquotesingle\textquotesingle CDF\textquotesingle\textquotesingle\ in Mathematica's SmoothHistogram). This last graph appears as a dashed red line, bounded between the cyan solid lines.

\bigskip

\bigskip

%\newpage

\section{Proofs}\label{Proofs}
%PPPPPPPPRRRRRRRRRROOOOOOOOOOOOOOFFFFFFFFFFFSSSSSS
\begin{lemma}\label{lemma1}
Let  $0\le x_1,\ldots,x_n <1$ and let $x_{\max}
=\max_{1\le k\le n}x_k$. 
If $x_{\max}\xrightarrow[n\to\infty]{} 0$ 
then, as~$n\to\infty$,
\begin{equation}%\label{lemma1.2}
  \prod_{k=1}^{n}\left(1-{x_k}\right)=\exp\left(-\sum_{k=1}^{n} x_k\right)
  +O\left( nx^2_{\max}\right).%A comment about the error: (i) It is $O_{-}$. (ii) We could have taken it to be O( \sum_{i=1}^n x^2_i). (iii) we could calculate it to be ...e...
\end{equation}
\end{lemma}
\begin{Proof}{}
%Recall  that \begin{equation}\label{3}   1+x\le e^x,\qquad x\in\bf{R}.\end{equation}
%$x< -1$, the \LHS of (\ref{3}) is negative, while the \RHS is positive. For $x\in[-1,0]$, consider the function $1+x-e^x$: At $-1$  its negative,  it  zeroes at $0$, and its  derivative, $1-e^x$ is non negative in this domain, which implies  it doesn't change sign. For $x\ge 0$:$$e^x=1+x+\frac{x^2}{2!}+\frac{x^3}{3!}+\cdots\ge 1+x.$$Altogether, (\ref{3}) is implied.
%By \citep[Sec 4.2, p.103]{Polya1952} we also have \begin{equation}    \label{4}e^{-x}\le 1-x+x^2,\qquad x\ge 0.\end{equation} %Consider the function $e^{-x}-1+x-x^2$. At $0$ it equals $0$. Now, consider its derivative: $-e^{-x}+1-2x$. By (\ref{3}) and for $x\ge 0$$$e^{-x}\ge 1-x\ge 1-2x.$$Thus, the function is non-increasing for $x\ge 0$, which implies (\ref{4}).
By 
\citep[p.103, (4.2.2)]{Polya1952} we have 
\begin{equation}\label{4a} 
e^{-x}= 1-x+\tfrac12x^2e^{-x\theta},\qquad x\ge 0,
\end{equation}where $0<\theta <1$. 
By (\ref{4a}),  for some $0< \theta_k <1$, $1\le k\le n$:
\begin{equation*}
    %\label{1}
\begin{aligned}
     1-x_k
    =e^{-x_k}-\tfrac12{x^2_k}e^{-\theta_k x_k}
    =e^{-x_k}\left(1-\tfrac12{x^2_k}e^{(1-\theta_k) x_k}\right).
    %%%%%%%%%%%%%\\
%    \ge & e^{-x_k}\left(1-\tfrac{e}{2}{x^2_k}\right)\ge e^{-x_k}\left(1-\tfrac{e}{2}{x^2_{\max}}\right).
\end{aligned}
\end{equation*}
Thus, %as $x_{\max}\to 0$
\begin{equation}\label{5}
\begin{aligned}
\prod_{k=1}^{n}\left(1-{x_k}\right)
&=\prod_{k=1}^{n}\left(e^{-x_{k}}\cdot\left(1-\tfrac12{x^2_k}e^{(1-\theta_k) x_k}\right)\right)\\
&=\exp\left(-\sum_{k=1}^{n} x_k\right)\prod_{k=1}^{n}\left(1-\tfrac12{x^2_k}e^{(1-\theta_k) x_k}\right).
\end{aligned}
\end{equation}
Consider the product on the \RHS of (\ref{5}). By Bernoulli's inequality:
\begin{equation}\label{6a}
\begin{aligned}
1\ge & \prod_{k=1}^{n}\left(1-\tfrac12{x^2_k}e^{(1-\theta_k) x_k}\right)
\ge \prod_{k=1}^{n}\left(1-\tfrac{e}{2}{x^2_{\max}}\right)
\ge  1-\tfrac{e}{2}\cdot n x^2_{\max} = 1+O\left(n x^2_{\max}\right).
%1&\ge\left(1-\tfrac{e}{2}{x^2_{\max}}\right)^n
%=\sum_{j=0}^n(-1)^j\binom{n}{j}\left(\frac{e}{2}\right)^j{x^{2j}_{\max}}\\
%=&1-n\tfrac{e}{2}{x^2_{\max}}+\sum_{j=1}^{n/2}\left(\binom{n}{2j}\left(\frac{e}{2}\right)^{2j}{x^{4j}_{\max}}
%-\binom{n}{2j+1}\left(\frac{e}{2}\right)^{2j+1}{x^{4j+2}_{\max}}\right)\\
%\ge & 1-n\tfrac{e}{2}x^{2}_{\max} =1+O\left(nx^{2}_{\max} \right).
\end{aligned}
\end{equation}
The lemma follows from (\ref{5}) and (\ref{6a}).

\qedForProof
\end{Proof}

\bigskip

\begin{lemma} \label{lemma2}
For fixed $a,\ll>0$ and  $c=O(\log\log n)$:
\begin{description}
\Item[{ $a.$}]
\begin{equation}\label{lemma2.a}
\begin{aligned}
&\left({\ec\log n}/{n}\right)^{\el\left(1-{\ll}/{n}\right)^{n-an/\log n}}\\
&\hspace{10mm}=\frac{\log n}{n}\cdot e^{-c-\la{}}
\left(1+{\ll a \log\log n}/{\log n}+\errLognly\right).
\end{aligned}
\end{equation}
\smallskip
\Item[{ $b.$}]
\begin{equation}\label{lemma2.new.a}
\begin{aligned}
\left({\ec\log n}/{n}\right)^{\tfrac{\ll}{n}\el\left(1-{\ll}/{n}\right)^{n-an/\log n}}
&=1-\frac{\ll}{n}\left(\log n -\log\log n+c+O(1)\right).
\end{aligned}
\end{equation}
\smallskip\Item[{ $c.$}]
\begin{equation}\label{lemma2.c}
\begin{aligned}
&\left(\ec\log n/n\right)^{\el\ql^{n-an \lnln/\log n}}\\
&\hspace{1cm}=\exp\left(-\left(\log n+(\la-1)\log\log n+c+O(\log^2\log n/\log n   )\right)\right).
\end{aligned}
\end{equation}
\end{description}
\end{lemma}
%\begin{remark}\label{remark1}In fact, we will use (\ref{lemma2.c}) with the exponent on the \LHS multiplied by $g(n)/n$ for the two functions $g(n)={\ll}\left({an(\log\log n -1)}/{\log n}+1\right)$ (in (\ref{57.1}) below) and for $g(n)=\lambda$ (in (\ref{53.2a}) below). In both cases, the argument of the $\exp$ on the \RHS is simply multiplied by $g(n)/n$ as well. \end{remark}
%\newpage
\begin{Proof}{\ of Lemma \ref{lemma2}}
\begin{description}
%aaaaaaaaaaaaaaaaaaaaaa
\item[{ $a.$}]
We start with the logarithm of the exponent on the \LHS of (\ref{lemma2.a}):
\begin{equation}\label{9}
\begin{aligned}
\log\left({\el\ql^{\an}}\right)
&=   \ll + {(\an)\log\ql}\\
&=   \ll+ {n(1-a/\log n)\left(-{\ll}/{n}+O(1/n^2)\right)}\\
&=   {\ll a/{\log n}+O(1/n)}.
\end{aligned}
\end{equation}
Hence:
\begin{equation}\label{10}
\begin{aligned}
\el\big(1-&{\ll}/{n}\big)^{n-an/\log n}\\
&=e^{\ll a/{\log n}+O(1/n)}\\
&=1+{\ll a}/{\log n}+O\left({1}/{n}\right)
+\tfrac12\left({\ll a}/{\log n}+O\left({1}/{n}\right)\right)^2+O\left({1}/{\log^3 n}\right)\\
&=1+{\ll a}/{\log n}+O\left({1}/{\log^2 n}\right)
.
\end{aligned}
\end{equation}
We shall deal separately with the factors $e^{-c}$ and $\log n/n$ in the base of the exponent on the \LHS of (\ref{lemma2.a}). For $e^{-c}$, by (\ref{10}): 
\begin{equation}\label{11}
\begin{aligned}
\left(e^{-c}
\right)^{\el\ql^{\an}}
&=\exp\left(-c\cdot\left(1+{\ll a}/{\log n}+O\left({1}/{\log^2 n}\right)\right)\right)\\
&=e^{-c}\exp\left(O\left({1}/{\log n}\right)\right)\\
&=e^{-c}\left(1+O\left({1}/{\log n}\right)\right).
\end{aligned}
\end{equation}
For the second factor $\log n/n$:
\begin{equation}\label{12}
\begin{aligned}
\left({\log n}/{n}
\right)^{\el\ql^{\an}}
&=\left({\log n}/{n}
\right)^{1+{\ll a}/{\log n}+O\left({1}/{\log^2 n}\right)}.
\end{aligned}
\end{equation}
Let us consider the logarithm of the \RHS of (\ref{12}):
\begin{equation}\label{12a}
\begin{aligned}
\log\Big(\big({\log n}&/{n}\big)^{1+{\ll a}/{\log n}+O\left({1}/{\log^2 n}\right)}\Big)\\
=\:&\left({1+{\ll a}/{\log n}+O\left({1}/{\log^2 n}\right)}\right)\left(-\log n+\log\log n \right)\\
=\:&-\log n+\log \log n-\la{}+{\ll a \log\log n}/{\log n}+\errLognly.
\end{aligned}
\end{equation}
Thus, by \Eq{12} and (\ref{12a}):
\begin{equation}\label{13}
\begin{aligned}
&\left({\log n}/{n}
\right)^{\el\ql^{\an}}\\
&\hspace{10mm}=\exp\left(-\log n +\log\log n-\la{}+\tfrac{\ll a \log\log n}{\log n}+\errLognly\right)\\
&\hspace{10mm}=\frac{\log n}{n}\cdot e^{-\la{}}\cdot\exp\left({\ll a \log\log n}/{\log n}+\errLognly\right)\\
&\hspace{10mm}=\frac{\log n}{n}\cdot e^{-\la{}}\cdot\left(1+{\ll a \log\log n}/{\log n}+\errLognly\right).
\end{aligned}
\end{equation}
By \Eq{11} and \Eq{13}:
\begin{equation*}%\label{lemma2.14}
\begin{aligned}
&\big(\ec\log
\,n/{n}\big)^{
    \el\left(1-{\ll}/{n}\right)^{n-an/\log n}}\\
&\hspace{1cm}=e^{-c}\left(1+\errLognly\right)\\
&\hspace{15mm}\cdot\frac{\log n}{n}\cdot e^{-\la{}}\left(1+{\ll a \log\log n}/{\log n}+\errLognly\right)\\
&\hspace{1cm}=\frac{\log n}{n}\cdot e^{-c-\la{}}\left(1+{\ll a \log\log n}/{\log n}+\errLognly\right).
\end{aligned}
\end{equation*}
\item[{ $b.$}]
%By (\ref{9}), similarly to (\ref{10}):  NOT SURE THIS IS NEEDED. I MIGHT JUST QUOTE (\ref{10}) IN THE NEXT SENTENCE.
%\begin{equation}\label{new.10}\begin{aligned}\frac{\ll}{n}\el\big(1-&{\ll}/{n}\big)^{n-an/\log n}%\\
%%&=\frac{\ll}{n}e^{\ll a/{\log n}+O(1/n)}\\&=\frac{\ll}{n}\left(1+{\ll a}/{\log n}+O\left({1}/{n}\right)+\tfrac12\left({\ll a}/{\log n}+O\left({1}/{n}\right)\right)^2+O\left({1}/{\log^3 n}\right)\right)\\
%&=\frac{\ll}{n}\left(1+{\ll a}/{\log n}+O\left({1}/{\log^2 n}\right)\right).\end{aligned}\end{equation}
Similarly to the proof of the previous part, we shall deal separately with the factors~$e^{-c}$ and $\log n/n$ in the base of the exponent on the \LHS of (\ref{lemma2.new.a}). For $e^{-c}$, by  (\ref{10})
\begin{equation*}%\label{new.11}
\begin{aligned}
\left(e^{-c}
\right)^{\tfrac{\ll}{n}\el\ql^{\an}}
&=\exp\left(-c\cdot\frac{\ll}{n}\left(1+{\ll a}/{\log n}+O\left({1}/{\log^2 n}\right)\right)\right)\\
&=1-c\cdot\frac{\ll}{n}\left(1+{\ll a}/{\log n}+O\left({1}/{\log^2 n}\right)\right)\\
&=1-c\ll/n+O\left({\lnln}/{(n\log n)}\right).
\end{aligned}
\end{equation*}
For the second factor $\log n/n$, by (\ref{12}), (\ref{12a}), and the first equality in (\ref{13}) 
%\begin{equation*}%\label{new.12a}\begin{aligned}&\log\Big(\big({\log n}/{n}\big)^{\tfrac{\ll}{n}\el\ql^{\an}}\Big)\\&\hspace{10mm}=\log\Big(\big({\log n}/{n}\big)^{\tfrac{\ll}{n}\left(1+{\ll a}/{\log n}+O\left({1}/{\log^2 n}\right)\right)}\Big)\\&\hspace{10mm}=\frac{\ll}{n}\left(-\log n+\log \log n-\la{}+{\ll a \log\log n}/{\log n}+\errLognly\right).\end{aligned}\end{equation*}Thus, 
\begin{equation*}%\label{new.13}
\begin{aligned}
&\left({\log n}/{n}
\right)^{\tfrac{\ll}{n}\el\ql^{\an}}\\
&\hspace{10mm}=\left(\left({\log n}/{n}
\right)^{\el\ql^{\an}}\right)^{\tfrac{\ll}{n}}\\
&\hspace{10mm}=\exp\left(\frac{\ll}{n}\left(-\log n +\log\log n-\la{}+\tfrac{\ll a \log\log n}{\log n}+\errLognly\right)\right)\\
&\hspace{10mm}=1+\frac{\ll}{n}\left(-\log n +\log\log n+O(1)\right).
\end{aligned}
\end{equation*}
Thus,%By \Eq{new.11} and \Eq{new.13}:
\begin{equation*}%\label{lemma2.14}
\begin{aligned}
&\big(\ec\log
\,n/{n}\big)^{\tfrac{\ll}{n}
    \el\left(1-{\ll}/{n}\right)^{n-an/\log n}}\\
&\hspace{1cm}=\left(1-\frac{c\ll}{n}+O\left({\lnln}/{(n\log n)}\right)\right)
\left(1-\frac{\ll}{n}\left(\log n -\log\log n+O(1)\right)\right)\\
&\hspace{1cm}= 1-\frac{\ll}{n}\left(\log n -\log\log n+c+O(1)\right).
\end{aligned}
\end{equation*}
%Going over the proof, one can check that, when changing the \LHS of (\ref{lemma2.new.a}) to  
%\begin{equation*}   \left({\ec\log n}/{n}\right)^{\tfrac{\ll}{n}\el\left(1-{\ll}/{n}\right)^{n-an/\log n-1}}\end{equation*}the \RHS does not change.

%bbbbbbbbbbbbbbbbbbb
\item[{ $c.$}]
Again we start with the logarithm of the exponent on the \LHS of (\ref{lemma2.c}):
\begin{equation*}%\label{9d}
\begin{aligned}
\log\left({\el\ql^{\anlln}}\right)
&=   \ll + {(\anlln)\log\ql}\\
&=   \ll+ {n(1-a\lnln/\log n)\left(-{\ll}/{n}+O(1/n^2)\right)}\\
&=   {\ll a\lnln/{\log n}+O(1/n)}.
\end{aligned}
\end{equation*}
Hence:
\begin{equation}\label{10d}
\begin{aligned}
&\el\left(1-{\ll}/{n}\right)^{\anlln}\\
&\hspace{10mm}= e^{\la \lnln/{\log n}+O(1/n)}\\
&\hspace{10mm}=\left(1+{\la \lnln}/{\log n}+O\left({1}/{n}\right)+O\left({\lln2}/{\log^2 n}\right)\right)\\
&\hspace{10mm}=\left(1+{\la \lnln}/{\log n}+O\left({\lln2}/{\log^2 n}\right)\right)
.
\end{aligned}
\end{equation}
Therefore,
\begin{equation}\label{54.3d} 
\begin{aligned}
&\left(\ec\right)^{\el\ql^{\anlln}}=\exp\left(-c+O(\lln2/\log n)\right).
\end{aligned}
\end{equation}
Also by (\ref{10d}),
\begin{equation*}
\begin{aligned}
&\log\left(\left(\log n/n\right)^{\el\ql^{\anlln}}\right)\\
&\hspace{1cm}=\left(1+\la \lnln/{\log n}+O(\lln2/\log^2n)\right)\left(\log\log n-\log n\right)\\
&\hspace{10mm}=-\left(\log n+\la \lnln-\log\log n
+O\left({\lln2}/{\log n}\right)\right).
\end{aligned}
\end{equation*}
Thus,
\begin{equation}\label{55.1d} 
\begin{aligned}
&\left(\log n/n\right)^{\el\ql^{\anlln}}\\
&\hspace{10mm}=\exp\left(-\left(\log n+(\la-1)\log\log n+O\left({\lln2 }/{\log n}\right)\right)\right).
\end{aligned}
\end{equation}
By  (\ref{54.3d}) and (\ref{55.1d}):
\begin{equation*}%\label{lemma2.d1} 
\begin{aligned}
&\left(\right.\ec\log n/n\left.\right)^{\el\ql^{\anlln}}\\
&\hspace{1cm}=\exp\left(- c+O(\lln2/\log n)\right)\\
&\hspace{1.5cm}\cdot\exp\left(-\left(\log n+(\la-1)\log\log n+O\left({\lln2 }/{\log n}\right)\right)\right)\\
&\hspace{1cm}=\exp\left(-\left(\log n+(\la-1)\log\log n+c+O(\log^2\log n/\log n   )\right)\right).
\end{aligned}
\end{equation*}
\end{description}
\qedForProof
\end{Proof}

\bigskip

\begin{Proof}{\ of Theorem \ref{ConvInDist'}}
It will be more convenient to work first with $T''=T'\cdot\ll/\el$, and then go back to  $T'$. 
We have:
\begin{equation}\label{T''}
\begin{aligned}
    {T''}=\frac{T-(\el/\ll)n(\log n-\log \log n)}{(\el/\ll)n}.
    \end{aligned}
\end{equation}
Denote by $F_{T''}$ its distribution function. For $c\in\bf{R}$:
\begin{equation}\label{1a}    \begin{aligned}
    F_{T''}(c)
    &=P\left(\frac{T-(\el/\ll)n(\log n-\log \log n)}{(\el/\ll)n}\le c\right)\\
    &=F_T\left(\tfrac{e^{\lambda}n}{\lambda} (\log n-\log \log n+c)\right).
    \end{aligned}
\end{equation}
By independence, for $t\ge 0$
\begin{equation}\label{1}
    \begin{aligned}
    \sizedSubLetter{F}{T}
    (t)&=P\left(\max_{1\le i\le n}\sizedSubLetter{T}{i}\le t \right)=\prod_{i=1}^{n}P\left(\sizedSubLetter{T}{i}\le t\right)\\
    &=\left(1-e^{-p_1t}\right)\left(1-e^{-p_2t}\right)\cdot\cdots\cdot\left(1-e^{-p_nt}\right)\\
    &=\prod_{i=1}^{n}\left(1-\exp\left({-(\lambda/n)(1-\lambda/n)^{i-1}t}\right)\right)\\
    &=\prod_{i=0}^{n-1}\left(1-\exp\left({-(\lambda/n)(1-\lambda/n)^it}\right)\right).
    \end{aligned}
\end{equation}

Let $t=(\el/\ll)n (\log n-\log \log n+c)$, where   $-
\aA\log\log\log n\le c\le \aB\log\log n$. (The need of dealing with unbounded values of $c$ arises from the proof of Theorem \ref{Expectation} below.) We start by estimating a typical term in the product on the \RHS of \Eq{1}.  First consider the exponent:
\begin{equation*}%\label{2a}
    \begin{aligned}
    -\tfrac{\ll}{n}\left(1-\tfrac{\ll}{n}\right)^it
    &=-\tfrac{\ll}{n}\left(1-\tfrac{\ll}{n}\right)^i\cdot{\tfrac{\el n}{\ll}(\log n-\log \log n+c)}\\
    &={(-\log n+\log \log n-c)}
    \cdot\el\left(1-{\ll}/{n}\right)^i.
    \end{aligned}
\end{equation*}
Thus,
\begin{equation}\label{2}
    \begin{aligned}
    1-\exp\left({-\tfrac{\ll}{n}\left(1-\tfrac{\ll}{n}\right)^it}\right)
    &=1-\exp\left({{(-\log n+\log \log n-c)}
    \cdot\el\left(1-{\ll}/{n}\right)^i}\right)\\
    &=1-\left({\ec\log n}/{n}\right)^{
    \el\left(1-{\ll}/{n}\right)^i}.
    \end{aligned}
\end{equation}
By \Eq{1}, \Eq{2}, and Lemma \ref{lemma1}
\begin{equation}\label{7}
\begin{aligned}
F_T\left(t\right)
=&\prod_{i=0}^{n-1}\left(1-\left({\ec\log n}/{n}\right)^{
    \el\ql^i}\right)\\
=&\exp\left({-
\sum_{i=0}^{n-1}
\left({\ec\log n}/{n}\right)^{    \el\ql^i}}\right)\\
    &
    +O\left(n\left({\ec\log n}/{n}\right)^{
    2\el\ql^{n-1}}\right).
\end{aligned}
\end{equation}

%WE WANT TO UNDERSTAND THE EXPONENT IN THE ERROR. NOW:
%$$\lambda+(n-1)\log(1-\lambda/n)=\\
%\lambda+(n-1)(-\lambda/n -\lambda^2/2n^2+O(1/n^3))\\
%=\ll/n-(n-1)\ll^2/2n^2+O()>-\ll^2/n$$

Consider the $i$-th addend in the sum in the exponent on the \RHS of (\ref{7}). Fix an $a\ge 2/\ll$, and let $i=n-an/\log n$. 
By Lemma \ref{lemma2}.a:
\begin{equation*}
\begin{aligned}
&\left({\ec\log n}/{n}\right)^{\el\left(1-{\ll}/{n}\right)^{i}}=\frac{\log n}{n}\cdot e^{-c-\la{}}\left(1+{\ll a \log\log n}/{\log n}+\errLognly\right).
\end{aligned}
\end{equation*}
Now, we estimate the whole sum in the exponent in the first addend on the \RHS of \Eq{7}. We split the sum into two. The first  consists of most of the addends of the sum, but, for large $a$, they contribute very little. 
The second  consists of the remaining minority, which accounts for most of the sum.
\begin{equation}\label{15}
\begin{aligned}
{
\sum_{i=0}^{n-1}
\left({\ec\log n}/{n}\right)^{    \el\ql^i}}
=&\:\sum_{0\le i\le n-{an}/{\log n}}
\left({\ec\log n}/{n}\right)^{    \el\ql^i}\\
&\:+\sum_{n-{an}/{\log n}<i\le n-1}% here we assume the reader understands that we shall not account for i twice. 
\left({\ec\log n}/{n}\right)^{    \el\ql^i}.
\end{aligned}
\end{equation}
We will  bound  the first sum on the \RHS
 of \Eq{15} both from above and from below. Let us start with an upper bound. Denote
 \begin{equation}\label{46.0}
 r_j=\left({\ec\log n}/{n}\right)^{
    \el\ql^{n-an/\log n -j}},\qquad 0\le j\le n-{an}/{\log n},    
 \end{equation}
so that:
\begin{equation}\label{46.1}
\begin{aligned}
\sum_{i=0}^{n-{an}/{\log n}}
\left({\ec\log n}/{n}\right)^{    \el\ql^i}
%&=\sum_{j=0}^{n-{an}/{\log n}}\left({\ec\log n}/{n}\right)^{\el\ql^{n-an/\log n -j}}\\
&=\sum_{j=0}^{n-{an}/{\log n}}
r_j.
\end{aligned}
\end{equation}
Clearly, $r_{j+1}\le r_j$ for ${0}\le j\le {n-an/\log n-1}$. Denote:
$$q_j=r_{j+1}/r_j,\qquad 0\le j\le n-{an}/{\log n}-1.
 $$Thus,
\begin{equation}\label{46.2a}
\begin{aligned}
q_j
&=\frac{\left({\ec\log n}/{n}\right)^{
    \el\left(1-{\ll}/{n}\right)^{n-an/\log n -j-1}}}{\left({\ec\log n}/{n}\right)^{
    \el\left(1-{\ll}/{n}\right)^{n-an/\log n -j}}}
=\left({\ec\log n}/{n}\right)^{g_1(\ll,a,n)},
\end{aligned}
\end{equation}
where
\begin{equation}\label{46.2b}
\begin{aligned}
{g_1(\ll,a,n)}
&={
    \el\ql^{n-an/\log n -j-1}}
    -{
    \el\ql^{n-an/\log n -j}}\\
%&={    \el\ql^{n-an/\log n -j-1}}\left(1-\ql\right)\\
&=\frac{\ll}{n}\cdot{
    \el\ql^{n-an/\log n -j-1}}.
\end{aligned}
\end{equation}
Therefore, by \Eq{46.2a} and \Eq{46.2b}:
\begin{equation}\label{46.2}
\begin{aligned}
q_j
&=\left({\ec\log n}/{n}\right)^{\tfrac{\ll}{n}\cdot{
    \el\ql^{n-an/\log n -j-1}}}=r_{j+1}^{\ll/n}.
\end{aligned}
\end{equation}
This implies that $q_{j+1}\le q_j$ for ${0}\le j\le {n-an/\log n-2}$. Hence:
\begin{equation*}%\label{46.3}
\begin{aligned}
\max_{0\le j\le n-{an}/{\log n-1}}q_j
&=q_0.
%=\left({\ec\log n}/{n}\right)^{\tfrac{\ll}{n}\cdot{\el\ql^{n-an/\log n-1}}}=\left(r_{1}\right)^{\ll/n}.
\end{aligned}
\end{equation*}
Thus,
\begin{equation}\label{47.1}
\begin{aligned}
\sum_{j=0}^{n-{an}/{\log n}}
r_j&\le\sum_{j=0}^{n-{an}/{\log n}}
r_0\cdot q_0^j\le\frac{r_0}{1-q_0}.
\end{aligned}
\end{equation}
We need to estimate the \RHS of (\ref{47.1}). By  Lemma \ref{lemma2}.a:
\begin{equation}\label{47.2}
\begin{multlined}
r_0=\frac{\log n}{n}\cdot e^{-c-\la{}}
\left(1+{\ll a \log\log n}/{\log n}+\errLognly{}\right).
\end{multlined}
\end{equation}
Now, for $q_0$, by (\ref{46.2}) and Lemma \ref{lemma2}.b (with the exponent slightly changed):
\begin{equation}\label{q0}
\begin{aligned}
q_0=1-\frac{\ll}{n}\left(\log n -\log\log n+c+O(1)\right).
\end{aligned}
\end{equation}
It follows that:
\begin{equation}\label{49.3}
\begin{aligned}
{1}/\big(1&-q_0\big)\\
&=\left(\frac{\ll}{n}\left(\log n-\log \log n +c +\errOne{}\right)\right)^{-1}\\
&=\frac{n}{\ll \log n}\left(1-\log \log n/\log n+c/\log n+ \errLognly\right)^{-1}\\
&=\frac{n}{\ll \log n}\left(1+\log \log n/\log n -c/\log n+\errLognly +O\left({\log^2\log n}/{\log^2 n}\right)\right)\\
&=\frac{n}{\ll \log n}\left(1+\log \log n/\log n -c/\log n+\errLognly\right).
\end{aligned}
\end{equation}
By (\ref{46.1}) and (\ref{47.1})-(\ref{49.3}), we obtain the following upper bound on the first sum on the \RHS of (\ref{15}):
\begin{equation}\label{49.4}
\begin{aligned}
{r_0}/\big(1-q_0\big)
=&\:\frac{\log n}{n}\cdot e^{-c-\la{}}\left(1+{\ll a \log\log n}/{\log n}+\errLognly\right)\\
&\cdot\frac{n}{\ll \log n}\left(1+\log \log n/\log n -c/\log n+\errLognly\right)\\
=&\:\frac{e^{-c-\la{}}}{\ll}\left(1+(\la+1){\log\log n}/{\log n}-c/\log n+\errLognly\right).
\end{aligned}
\end{equation}
We now establish a lower bound on the first sum on the \RHS of (\ref{15}). By (\ref{46.0}) and the change of variable $j=\an-i$,
\begin{equation}\label{53.1a}
\begin{aligned}
\sum_{i=0}^{n-{an}/{\log n}}
\left({\ec\log n}/{n}\right)^{    \el\ql^i}
&\ge \sum_{i=n-{an}\log\log n/{\log n}-1}^{\an}\left({\ec\log n}/{n}\right)^{    \el\ql^i}\\
&= \sum_{i=n-{an}\log\log n/{\log n}-1}^{\an}r_{\an-i}
=\sum_{j=0}^{{an}(\log\log n-1)/{\log n}}
r_j\\
&=r_0+\sum_{j=0}^{{an}(\log\log n-1)/{\log n}-1}
r_0\cdot q_0\cdots q_j.
\end{aligned}
\end{equation}
By (\ref{53.1a}), and since $q_j$ decreases as a function of $j$,
\begin{equation}\label{53.1}
\begin{aligned}
\sum_{i=0}^{n-{an}/{\log n}}
\left({\ec\log n}/{n}\right)^{    \el\ql^i}
&\ge\sum_{j=0}^{{an}(\log\log n-1)/{\log n}}
r_0\left(q_{{an}(\log\log n-1)/{\log n}-1}\right)^j\\
&=r_0\cdot\frac{1-\left(q_{{an}(\log\log n-1)/{\log n}-1}\right)^{{an}(\log\log n-1)/{\log n}+1}}{1-q_{{an}(\log\log n-1)/{\log n}-1}}.
\end{aligned}
\end{equation}
We will start with the second addend in the numerator on the \RHS of (\ref{53.1}). By  Lemma \ref{lemma2}.c,  and as $\la\ge2$:
\begin{equation}\label{57.1}
\begin{aligned}
&\left(q_{\anllnq-1}\right)^{\anllnq+1}\\
&\hspace{1cm}=\left(\ec\log n/n\right)^{\tfrac{\ll}{n}(\anllnq+1)\el\ql^\anlln}\\
&\hspace{1cm}=\exp\left(-\frac{\ll}{n}\left(\tfrac{an(\log\log n -1)}{\log n}+1\right)
\left(\log n+(\la-1)\log\log n+c+O\left(\tfrac{c\lnln}{\log n}\right)\right)\right)\\
&\hspace{1cm}=\exp\left(-\la(\lnln-1)+O\left({\lln2}/{\log n }\right)\right)\\
&\hspace{1cm}=\left({e}/{\log n}\right)^{\la}\left(1+O\left({\lln2/\log n}\right)\right)
=%\left({e}/{\log n}\right)^{\la}+
O\left({1/\log^2 n }\right).
\end{aligned}
\end{equation}
For the denominator on the \RHS of (\ref{53.1}), by (\ref{46.2}) and Lemma \ref{lemma2}.c: 
\begin{equation}\label{53.2a}
    \begin{aligned}
    q_{{\anllnq}-1}
    &=\left(\ec\log n/n\right)^{\tfrac{\ll}{n}\el\ql^\anlln}\\
    &=\exp\left(-\frac{\ll}{n}\left(\log n+(\la-1)\log\log n+c+O\left(\tfrac{c\log\log n}{\log n}\right)\right)\right)\\
    &=1-\frac{\ll}{n}\left(\log n+(\la-1)\log\log n+c+O\left(\tfrac{c\log\log n}{\log n}\right)\right).
    \end{aligned}
\end{equation}
Thus,
\begin{equation}\label{57.2}
\begin{aligned}
&{1}/\left({1-q_{{\anllnq}-1}}\right)\\
&\hspace{1cm}=\left(\frac{\ll}{n}\left(\log n+(\la-1)\log\log n+c+O\left(\tfrac{c\log\log n}{\log n}\right)\right)\right)^{-1}\\
&\hspace{1cm}=\frac{n}{\ll \log n}\left(1+\left(\la-1\right)\lnln/\log n +c/\log n+O\left(\tfrac{c\log\log n}{\log^2 n}\right)\right)^{-1}\\
&\hspace{1cm}=\frac{n}{\ll \log n}\left(1-\left(\la-1\right)\lnln/\log n-c/\log n+O\left(\tfrac{\log^2\log n}{\log^2 n}\right)\right).
\end{aligned}
\end{equation}
By (\ref{47.2}), (\ref{57.1}), and (\ref{57.2}), the lower bound on the first sum on the \RHS of (\ref{15}) is
\begin{equation}\label{58.2}
\begin{aligned}
&{r_0}\left({1-\left(q_{\anllnq-1}\right)^{\anllnq+1}}\right)/\left({1-q_{\anllnq-1}}\right)\\
&\hspace{1cm}=\frac{\log n}{n}\cdot e^{-c-\la{}}\left(1+{\ll a \log\log n}/{\log n}+\errLognly\right)
\cdot\left(1+O\left({1}/{\log^2 n }\right)\right)\\
&\hspace{1.5cm}\cdot\frac{n}{\ll \log n}\left(1-\left(\la-1\right)\lnln/\log n -c/\log n+O\left(\tfrac{\log^2\log n}{\log^2 n}\right)\right)\\
&\hspace{1cm}=\frac{e^{-c-\la{}}}{\ll}\left(1+\lnln/\log n-c/\log n+\errLognly\right).
\end{aligned}
\end{equation}
Now let us consider the second sum on the \RHS of (\ref{15}). %Note that 
%\begin{equation*}\left({\ec\log n}/{n}\right)^{    \el\ql^x}\le \left({\ec\log n}/{n}\right)^{    \el\ql^{x+1}},\qquad x\in\bf{R}.\end{equation*}
Clearly,
\begin{equation}\label{17}
\begin{aligned}
\int\limits_{\integrallimitL +1}^{n-1}\integrandbase{}^{    \el\ql^x}dx
&\le \sum_{\integrallimitL<i\le n-1}
\integrandbase{}^{    \el\ql^i}\\
&\le \int\limits_{\integrallimitL}^{n}\integrandbase{}^{    \el\ql^x}dx.
\end{aligned}
\end{equation}
Furthermore, 
\begin{equation}\label{18}
\begin{aligned}
0&\le\int\limits_{\integrallimitL }^{\integrallimitL +1}\integrandbase{}^{    \el\ql^x}dx
\le\int\limits_{n-1}^{n}\integrandbase{}^{    \el\ql^x}dx\\
&\le \integrandbase{}^{    \el\ql^n}
\le \integrandbase{}^{1/2}\le O\left(\log n/\sqrt{n}\right).
\end{aligned}
\end{equation}
By (\ref{17}) and (\ref{18}): 
\begin{equation}\label{19}
\begin{aligned}
\sum_{\integrallimitL<i\le n-1}
\integrandbase{}^{    \el\ql^i}
=& \int\limits_{\integrallimitL}^{n}\integrandbase{}^{    \el\ql^x}dx\\ &\:+O\left(\log n/\sqrt{n}\right).
\end{aligned}
\end{equation}
Consider the integral on the \RHS of (\ref{19}). By the change of variables $$y=(1-x/n)\log n,\quad x=n-yn/\log n,\quad dx=-n/\log n\,dy,$$ we obtain
\begin{equation*}\label{19a}
\begin{aligned}
\int\limits_{\integrallimitL}^{n}\integrandbase{}^{    \el\ql^x}dx
&= -\frac{n}{\log n}\int\limits_{a}^{0}\integrandbase^{\el\ql^{n-yn/\log n}}dy\\
&= \frac{n}{\log n}\int\limits_{0}^{a}\integrandbase^{\el\ql^{n-yn/\log n}}dy.
\end{aligned}
\end{equation*}
By Lemma \ref{lemma2}.a
\begin{equation}\label{39.4}
\begin{aligned}
&\int\limits_{\integrallimitL}^{n}\integrandbase{}^{    \el\ql^x}dx\\
&\hspace{1cm}= \frac{n}{\log n}\int\limits_{0}^{a}\frac{\log n}{n}\cdot e^{-c-\ll y}
\left(1+{\ll y \log\log n}/{\log n}+\errLognly\right)dy\\
&\hspace{1cm}= e^{-c}\left(1+\errLognly\right)\int\limits_{0}^{a} e^{-\ll y}dy+ \frac{\ll e^{-c}\lnln}{\log n} \int\limits_{0}^{a}y e^{-\ll y}dy\\
&\hspace{1cm}= e^{-c}\left(1+\errLognly\right)\cdot \frac{1}{\ll}\left(1-e^{-\la}\right)
+ \frac{\ll e^{-c}\lnln}{\log n} \cdot \frac{1}{\ll^2}\left(1-e^{-\la}\left(1+\la\right)\right)\\
&\hspace{1cm}
=\frac{ e^{-c}}{\ll} \left(1-e^{-\la}+\frac{\log\log n}{\log n}\cdot \left(1-e^{-\la}\left(1+\la\right)\right)+\errLognly\right).
\end{aligned}
\end{equation}
By (\ref{15}), (\ref{58.2}) and (\ref{39.4}):%, by (\ref{15}), (\ref{46.1}), (\ref{47.1}), (\ref{49.4}), (\ref{19}), and (\ref{39.4}),
 \begin{equation}\label{69.1}
\begin{aligned}
&\sum_{i=0}^{n-1}\integrandbase{}^{    \el\ql^i}\\
&\hspace{1cm}\le\frac{e^{-c-\la{}}}{\ll}\left(1+(\la+1)\cdot\lglgFlgly-c/\log n+\errLognly\right)\\
&\hspace{1.5cm}
+\frac{ e^{-c}}{\ll} \left(1-e^{-\la}+\frac{\log\log n}{\log n}\cdot \left(1-e^{-\la}\left(1+\la\right)\right)+\errLognly\right)\\
&\hspace{1cm}
=\frac{ e^{-c}}{\ll} \left(1+{\lnln}/{\log n}-ce^{-\la}/\log n+\errLognly\right),
\end{aligned}
\end{equation}
For the lower bound we get:%, by (\ref{15}), (\ref{53.1}), (\ref{58.2}),  (\ref{19}), and (\ref{39.4}),
\begin{equation}\label{60.1}
\begin{aligned}
&\sum_{i=0}^{n-1}
\integrandbase{}^{    \el\ql^i}\\
&\hspace{1cm}\ge\frac{e^{-c-\la{}}}{\ll}\left(1+\lnln/\log n-c/\log n+\errLognly\right)\\
&\hspace{1.5cm}+\frac{ e^{-c}}{\ll} \left(1-e^{-\la}+\frac{\log\log n}{\log n}\cdot \left(1-e^{-\la}\left(1+\la\right)\right)+\errLognly\right)\\
&\hspace{1cm}=\frac{ e^{-c}}{\ll} \left(1+\left(1-\la e^{-\la{}}\right){\lnln}/{\log n}-ce^{-\la}/\log n+\errLognly\right).
\end{aligned}
\end{equation}
We now use the above bounds to obtain corresponding bounds on $F_{T''}$.  For the lower bound, by (\ref{1a}), (\ref{7}), (\ref{69.1}), 
and since $c\ge -2\log\log\log n$,
\begin{equation}\label{69.3}
\begin{aligned}
F_{T''}(c)&=
F_T\left(\frac{n\el}{\ll} (\log n-\log \log n+c)\right)\\
&\ge\exp\left(-\frac{ e^{-c}}{\ll} \left(1+\lglgFlgly-ce^{-\la}/\log n+\errLognly\right)\right)\\
&\hspace{10pt}+O\left(n\left({\ec\log n}/{n}\right)^{ 2\el\ql^{n-1}}\right)\\
&=e^{-e^{-c}/\ll}\cdot\exp\left(-\frac{ e^{-c}}{\ll} \left(\lglgFlgly-ce^{-\la}/\log n+\errLognly\right)\right)\\
&\hspace{10pt}+O\left(n\left({\lln2\log n}/{n}\right)^{ 2-o(1)}\right)\\
&=e^{-e^{-c}/\ll}\cdot\exp\left(-\frac{ e^{-c}}{\ll} \left(\lglgFlgly-ce^{-\la}/\log n+\errLognly\right)\right).
\end{aligned}
\end{equation}
Similarly, for the upper bound,  by (\ref{1a}), (\ref{7}), and (\ref{60.1}):
\begin{equation}\label{70.0} 
\begin{aligned}
F_{T''}(c)
&\le\exp\left(-\frac{ e^{-c}}{\ll} \left(1+(1-\la e^{-\la})\lglgFlgly-ce^{-\la}/\log +\errLognly\right)\right).
\end{aligned}
\end{equation}
Thus, by (\ref{69.3}) and (\ref{70.0}), for any fixed $c$,
\begin{equation}\label{F_T''}
\begin{aligned}
F_{T''}(c)
%=F_T\left(\frac{n\el}{\ll} (\log n-\log \log n+c)\right)
=e^{-e^{-c}/\ll}+O \left(\lnln/{\log n}\right).
\end{aligned}
\end{equation}
Consequently
\begin{equation}\label{F_T'label2}
\begin{aligned}
    {F}_{T'}(c)
    &=P\left(\frac{T-(\el/\ll)\cdot n(\log n-\log \log n)}{n}\le c\right)\\
   %&=P\left(\frac{T-(\el/\ll)\cdot n(\log n-\log \log n)}{(\el/\ll)n}\le c\ll/\el\right)\\
    &=P\left(T\le \frac{\el}{\ll}\cdot n\left(\log n-\log \log n+c\ll/\el\right)\right)\\
    &=F_{T''}(c\ll/\el)
    =\exp\left(-e^{-c\ll/\el}/\ll\right)+O\left({\log\log n}/{\log n}\right)\\
    &=\exp\left(-e^{-\left(c-\left(-\el \log\ll/\ll\right)\right)/\left(\el/\ll\right)}\right)+O\left({\log\log n}/{\log n}\right),\qquad c\in\bf{R}.
\end{aligned}
\end{equation}

\qedForProof
\end{Proof}

\bigskip

%llllllllllleeeeeeeeeemmmmmmmmmmmmaaaaaaaaaaa
%llllllllllleeeeeeeeeemmmmmmmmmmmmaaaaaaaaaaa
%llllllllllleeeeeeeeeemmmmmmmmmmmmaaaaaaaaaaa
The following lemma will be used in the proof of Theorem \ref{F_D}, %** 
and may be of independent interest. 
\begin{lemma}\label{lemmaF_DF_T}For sufficiently large $d$:
\begin{description}
\item[{ $a.$}] $F_D(d)\le F_T\left(d+d^{3/4}\right)+1/\left(d+d^{3/4}\right)^{1/3}.$
\item[{ $b.$}] $F_D(d)\ge F_T\left(d-d^{3/4}\right)-1/\left(d-d^{3/4}\right)^{1/3}.$
\end{description}
\end{lemma}

\begin{Proof}{}
Let us construct a coupling of $D$ and $T$.

Consider the process of coupon arrivals under the continuous model. We take into account  real coupons as well as dummy coupons. %(Of course, in principle we are interested only in the time until we get all the real coupons. However, as explained in Section \ref{introduction1} with respect to the discrete process, it is inconsequential whether we require the dummy coupon to have arrived or not.) 

We present the continuous model, discussed in Section \ref{sec:Results}, in a somewhat different way. Suppose we get coupons according to a Poisson process with rate $1$, where each coupon is of type $0$ with probability $p_0$, of type $1$ with probability $p_1$, and so on. It is readily seen that the process is equivalent to the one in Section \ref{sec:Results} (where now we add a flow for  dummy coupons, with inter-arrival times $T_0$ distributed $\Exp(p_0)$). Namely, each coupon $i$ by itself is obtained according to a Poisson process with rate $p_i$, and the processes are independent for the various $i$-s.
Let $\widetilde{T}$ be the time until  all real coupon types have been received, and  $\widetilde{D}$  the total number of coupons, both real and dummy, received in the process.
Clearly, $T$ and $\widetilde{T}$ have the same distribution, and the same applies to  $D$ and  $\widetilde{D}$.

Consider the number $\N(t)$ of coupons (real or dummy) arriving until time~$t$ in the continuous process. By \citep[Ch.7]{Ross2009}, the variable  $\N(t)$ is Poisson distributed with parameter~$t$. By Chebyshev's inequality, the probability that we receive less than $t-t^{2/3}$ coupons  until time~$t$ is bounded as follows:
\begin{equation}\label{110.1}
    \begin{aligned}
    P\left(\N\left(t\right)\le t-t^{2/3}\right)%&=P\left(\N\left(t\right)-t\le-t^{2/3}\right)\le P\left(\left|\N(t)-t\right|\ge t^{2/3}\right)\\
    &\le V\left(\N\left(t\right)\right)/\left(t^{2/3}\right)^2=t/t^{4/3}=1/t^{1/3}.
    \end{aligned}
\end{equation}
By (\ref{110.1}),
\begin{equation}\label{110.1a}
    \begin{aligned}
    F_{\widetilde{T}}(t)&= P\left(\N\left(t\right)
    \le t-t^{2/3},{\widetilde{T}}\le t\right)+P\left(\N\left(t\right)>t-t^{2/3},{\widetilde{T}\le t}\right)\\
    &\ge P\left(\N\left(t\right)>t-t^{2/3}\right)\cdot P\left({\widetilde{T}}\le t\vert \N\left(t\right)>t-t^{2/3}\right)\\
    &\ge \left(1-1/t^{1/3}\right)\cdot P\left({\widetilde{D}}\le t-t^{2/3}\right)\ge  F_{\widetilde{D}}\left(t-t^{2/3}\right)-1/t^{1/3}.
    \end{aligned}
\end{equation}
Let $t=d+d^{3/4}$. For sufficiently large $d$:
\begin{equation}\label{110.3}
    t-t^{2/3}=d+d^{3/4}-\left(d+d^{3/4}\right)^{2/3}\ge d+d^{3/4}-\left(2d\right)^{2/3}\ge d.
\end{equation}
By (\ref{110.1a}) and (\ref{110.3}),
\begin{equation*}%\label{110.2}
    F_D(d)=F_{\widetilde{D}}(d)\le F_{\widetilde{D}}\left(t-t^{2/3}\right)\le F_{\widetilde{T}}\left(t\right)+1/t^{1/3}=F_{{T}}\left(t\right)+1/t^{1/3},
\end{equation*}
which proves the first part of the lemma.

To bound $F_{\widetilde{D}}\left(d\right)$ from below, we proceed in a similar way. We have
\begin{equation*}%\label{105.8}
    \begin{aligned}
    P\left(\N\left(t\right)\ge t+t^{2/3}\right)\le 1/t^{1/3},
    \end{aligned}
\end{equation*}
and:
\begin{equation*}%\label{105.7}
    \begin{aligned}
    F_{\widetilde{T}}(t)&= P\left(\N\left(t\right)
    \ge t+t^{2/3},{\widetilde{T}}\le t\right)+P\left(\N\left(t\right)<t+t^{2/3},{\widetilde{T}}\le t\right)\\
    &\le P\left(\N\left(t\right)\ge t+t^{2/3}\right)
    +P\left(\N\left(t\right)< t+t^{2/3}\right)\cdot P\left({\widetilde{T}}\le t\vert\N\left(t\right)< t+t^{2/3}\right)\\
    &\le P\left(\N\left(t\right)\ge t+t^{2/3}\right)
    + P\left({\widetilde{T}}\le t\vert\N\left(t\right)< t+t^{2/3}\right)\\
    &\le P\left(\N\left(t\right)\ge t+t^{2/3}\right)
    + F_{\widetilde{D}}\left( t+t^{2/3}\right).
    \end{aligned}
\end{equation*}

\begin{comment}
Consider the first addend on the \RHS of (\ref{105.7}). Similarly to (\ref{110.1})
\begin{equation}\label{105.8}
    \begin{aligned}
    P\left(\N\left(t\right)\ge t+t^{2/3}\right)&=P\left(\N\left(t\right)-t\ge t^{2/3}\right)\le P\left(\left|\N(t)-t\right|\ge t^{2/3}\right)\\
    &\le V\left(\N\left(t\right)\right)/\left(t^{2/3}\right)^2=t/t^{4/3}=1/t^{1/3}.
    \end{aligned}
\end{equation}
For some large $d$, (similarly to (\ref{110.3})) take $t=d-d^{3/4}$. Thus,
\begin{equation}\label{109.11}
    t+t^{2/3}=d-d^{3/4}+\left(d-d^{3/4}\right)^{2/3}\le d-d^{3/4}+d^{2/3}\le d.
\end{equation}
By (\ref{105.7}) and (\ref{109.11}):
\begin{equation}\label{109.12}
    F_D(d)\ge F_D\left(t+t^{2/3}\right)\ge F_T\left(t\right)-1/t^{1/3}=F_T\left(d-d^{3/4}\right)+1/\left(d-d^{3/4}\right)^{1/3}.
\end{equation}
\end{comment}

Taking $t=d-d^{3/4}$, we easily complete the proof.

\qedForProof
\end{Proof}

\bigskip

%%%%%%%%%%%%%%%%%%%%%%%%%%%%%%
\begin{Proof}{ of Theorem \ref{F_D}}
We only prove (\ref{F_D'.2}), as it clearly  implies (\ref{F_D'.1}). It will be more convenient to work first with $D''=D'\cdot\ll/\el$, and then return to  $D'$. We have:
\begin{equation*}
\begin{aligned}
    {D''}=\frac{D-(\el/\ll)n(\log n-\log \log n)}{(\el/\ll)n}.
    \end{aligned}
\end{equation*}
Similarly to (\ref{1a}), the distribution functions $F_{D''}$ and $F_D$ are related by:
\begin{equation}\label{105.1a}    \begin{aligned}
    F_{D''}(d'')
    =F_D\left(\tfrac{ne^{\lambda}}{\lambda} (\log n-\log \log n+d'')\right), \qquad d''\in\R.
    \end{aligned}
\end{equation}
For large $n$, denote
\begin{equation}\label{105.1c}
    d_n=\tfrac{ne^{\lambda}}{\lambda} (\log n-\log \log n+d'').
\end{equation}
By Lemma \ref{lemmaF_DF_T}.a,
\begin{equation}\label{105.2}
    F_D(d_n)\le F_T\left(d_n+d_n^{3/4}\right)+1/\left(d_n+d_n^{3/4}\right)^{1/3}.
\end{equation}
Consider the first term on the \RHS of (\ref{105.2}):
\begin{equation}\label{105.2a}
    \begin{aligned}
        d_n+d_n^{3/4}&=\tfrac{ne^{\lambda}}{\lambda} \left(\log n-\log \log n+d''\right)+\left(\tfrac{ne^{\lambda}}{\lambda} \left(\log n-\log \log n+d''\right)\right)^{3/4}\\
        &=\tfrac{ne^{\lambda}}{\lambda} \left(\log n-\log \log n+d''+O\left(\log n/n^{1/4}\right)\right).
    \end{aligned}
\end{equation}
Consider the variable
\begin{equation*}
\begin{aligned}
    {T''}=\frac{T-(\el/\ll)n(\log n-\log \log n)}{(\el/\ll)n},
    \end{aligned}
\end{equation*}
defined in (\ref{T''}).
By (\ref{1a}) and (\ref{F_T''}), for fixed $t''\in\bf{R}$: 
\begin{equation}\label{105.1b}
\begin{aligned}
F_{T''}(t'')&=F_T\left(\tfrac{ne^{\lambda}}{\lambda} \left(\log n-\log \log n+t''\right)\right)\\
&=\exp\left({-e^{-t''}/\ll}\right)+O \left(\lnln/{\log n}\right).
\end{aligned}
\end{equation}
By (\ref{105.2a}) and (\ref{105.1b}),
\begin{equation}\label{105.3}
    \begin{aligned}
        F_T\left(d_n+d_n^{3/4}\right)&=F_T\left(\tfrac{ne^{\lambda}}{\lambda} \left(\log n-\log \log n+d''+O\left(\log n/n^{1/4}\right)\right)\right)\\
        &=\exp\left({-e^{-d''+O\left(\log n/n^{1/4}\right)}/\ll}\right)+O \left(\lnln/{\log n}\right).
    \end{aligned}
\end{equation}
Now
\begin{equation}\label{105.3a}
    \begin{aligned}
        e^{-d''+O\left(\log n/n^{1/4}\right)}&=e^{-d''}\cdot e^{O\left(\log n/n^{1/4}\right)}
        =e^{-d''}\cdot\left( 1+{O\left(\log n/n^{1/4}\right)}\right)\\    &=e^{-d''} +{O\left(\log n/n^{1/4}\right)},    \end{aligned}
\end{equation}
so that the first term on the right-hand side of (\ref{105.3}) is:
\begin{equation}\label{105.4}
    \begin{aligned}
     \exp\left({-e^{-d''+O\left(\log n/n^{1/4}\right)}/\ll}\right)
     &=\exp\left({-e^{-d''} /\ll+{O\left(\log n/n^{1/4}\right)}}\right)\\   
    &=\exp\left({-e^{-d''} /\ll}\right)\cdot\left(1+{O\left(\log n/n^{1/4}\right)}\right)\\   
    &=\exp\left({-e^{-d''} /\ll}\right)+{O\left(\log n/n^{1/4}\right)}.
    \end{aligned}
\end{equation}
By (\ref{105.3}) and (\ref{105.4}),
\begin{equation}\label{105.5}
    \begin{aligned}
        F_T\left(d_n+d_n^{3/4}\right)&=\exp\left({-e^{-d''}/\ll}\right)+O \left(\lnln/{\log n}\right).
    \end{aligned}
\end{equation}
\begin{comment}
By (\ref{105.2a}), the second term on the \RHS of (\ref{105.2}) is,
\begin{equation}\label{105.5a}
    \begin{aligned}
        1/\left(d_n+d_n^{3/4}\right)^{1/3}&=1/\left(\tfrac{ne^{\lambda}}{\lambda} \left(\log n-\log \log n+d''+O\left(\log n/n^{1/4}\right)\right)\right)^{1/3}\\
        &=O\left(\left(\log n/n\right)^{1/3}\right).
    \end{aligned}
\end{equation}
\end{comment}
Thus, by (\ref{105.1a})-%, (\ref{105.1c}) 
(\ref{105.2}) and (\ref{105.5}), %and (\ref{105.5a})
\begin{equation}\label{105.6}
    \begin{aligned}
        F_{D''}\left(d''\right)&=F_{D}\left(d_n\right)\le\exp\left({-e^{-d''}/\ll}\right)+O \left(\lnln/{\log n}\right).
    \end{aligned}
\end{equation}
Similarly, by Lemma \ref{lemmaF_DF_T}.b we get the analogue of (\ref{105.6}), with the inequality reversed. Altogether:
\begin{equation*}\label{F_D''}
\begin{aligned}
F_{D''}(d'')
=\exp\left({-e^{-d''}/\ll}\right)+O \left(\lnln/{\log n}\right).
\end{aligned}
\end{equation*}
Similarly to (\ref{F_T'label2}), we get our claim.

\qedForProof
\end{Proof}

\begin{lemma}\label{LemmaForExpectation}
Let $\ll>0$.
\begin{description}
\item[{ $a.$}]There exists  a $\thA>0$ such that, for   every $ c\le 0$, 
\begin{equation*}
    \begin{aligned}  \sum_{i=0}^{n-1}\left({\ec\log n}/{n}\right)^{    \el\ql^i}
\ge e^{-c(1-o(1))}\cdot \thA.
    \end{aligned}
\end{equation*}
\item[{ $b.$}]There exists  a $\thB>0$ such that, for every  $c\ge 0$ 
\begin{equation*}
    \begin{aligned}  \sum_{i=0}^{n-1}\left({\ec\log n}/{n}\right)^{    \el\ql^i}
\le e^{-c(1-o(1))}\cdot \thB\log\log n.
    \end{aligned}
\end{equation*}
\end{description}
\end{lemma}
\begin{Proof}{\ of Lemma \ref{LemmaForExpectation}}
\begin{description}
%aaaaaaaaaaaaaaaaaaaaaa
\item[{ $a.$}]
As %$-\log n +\lnln \le c\le -\log\log n$ we have 
$1  \le e^{-c}$,
\begin{equation}\label{200.0}
    \begin{aligned}  \sum_{i=0}^{n-1}\left({\ec\log n}/{n}\right)^{    \el\ql^i}
\ge e^{-c(1-o(1))}\sum_{i=0}^{n-1}\left({\log n}/{n}\right)^{    \el\ql^i}.
    \end{aligned}
\end{equation}
Consider the sum on the \RHS of (\ref{200.0}). For $g(n)=n-n/\log n$
\begin{equation}\label{200.1}
    \begin{aligned}  \sum_{i=0}^{n-1}\left(\frac{\log n}{n}\right)^{    \el\ql^i}&\ge
%    \sum_{0\le i< g(n)}\left(\frac{\log n}{n}\right)^{    \el\ql^i}+
    \sum_{g(n)\le i\le n-1}\left(\frac{\log n}{n}\right)^{    \el\ql^i}\\
    &\ge %g(n)\cdot\left(\log n/n\right)^{e^{\ll}}+
    \left(n-g(n)\right)\cdot\left(\log n/n\right)^{    \el\ql^{g(n)}}.
    \end{aligned}
\end{equation}
%For the first addend on the \RHS of (\ref{200.1}),
%\begin{equation}\label{200.4}   \begin{aligned}   g(n)\cdot\left(\frac{\log n}{n}\right)^{e^{\ll}}   &=\left(n-\frac{n}{\log n}\right)\cdot\left(\frac{\log n}{n}\right)^{e^{\ll}}=\left(1-\frac{1}{\log n}\right)\cdot\frac{\log^{e^{\ll}} n}{n^{e^{\ll}-1}}.\end{aligned}
%\end{equation}
We start with the logarithm of the exponent of the second factor on the \RHS of (\ref{200.1}):
\begin{equation*}%\label{}
    \begin{aligned}   
    \log \left(    \el\ql^{g(n)}\right)
    &=\ll+g(n)\log\ql\\
    &=\ll+\left(n-n/\log n\right)\left(-\ll/n+O(1/n^2)\right)\\
    &=\ll-\ll+\ll/\log n+O(1/n)\\
    &=\ll/\log n+O(1/n).
    \end{aligned}
\end{equation*}
Thus,
\begin{equation*}%\label{200.2}
    \begin{aligned}   
    \el\ql^{g(n)}
    &=\exp\left(\ll/\log n-\tfrac12\ll^2/n+O(1/(n\log n))\right)\\
    &=1+\ll/\log n+O(1/\log^2 n).
    \end{aligned}
\end{equation*}
Hence,
\begin{equation*}%\label{}
    \begin{aligned}   
    &\log \left(    \left(n-g(n)\right)\cdot\left(\log n/n\right)^{    \el\ql^{g(n)}}\right)\\
    &\hspace{10mm}=\log\left(n-g(n)\right)+\el\ql^{g(n)}\left(-\log n+\log\log n\right)\\
    &\hspace{10mm}=\log n-\lnln+\left(1+\ll/\log n+O(1/\log^2 n)\right)\left(-\log n+\log\log n\right)\\
    &\hspace{10mm}=-\ll+\ll\lnln/\log n+O(1/\log n).
    \end{aligned}
\end{equation*}
Thus,
\begin{equation}\label{200.3}
    \begin{aligned}   
    \left(n-g(n)\right)\cdot\left(\log n/n\right)^{    \el\ql^{g(n)}}
    &=\exp\left(-\ll+\ll\lnln/\log n+O(1/\log n)\right)\\
    &=e^{-\ll}\left(1+\ll\lnln/\log n+O(1/\log n)\right).
    \end{aligned}
\end{equation}
By (\ref{200.1}) and (\ref{200.3}):
\begin{equation}\label{200.5}
    \begin{aligned}  &\sum_{i=0}^{n-1}\left({\log n}/{n}\right)^{    \el\ql^i}
    %&\hspace{10mm}\ge\left(1-\frac{1}{\log n}\right)\cdot\frac{\log^{e^{\ll}} n}{n^{e^{\ll}-1}}+e^{-\ll}\left(1+\ll\lnln/\log n+O(1/\log n)\right)\\
    %&\hspace{10mm}\ge e^{-\ll}\\
    &\ge%\left(1-\frac{1}{\log n}\right)\cdot\frac{\log^{e^{\ll}} n}{n^{e^{\ll}-1}}+
    e^{-\ll}\left(1+\frac{\ll\lnln}{\log n} +O\left(\frac{1}{\log n}\right)\right)\ge e^{-\ll}.
    \end{aligned}
\end{equation}
The claim follows from (\ref{200.0}) and (\ref{200.5})  for $\thA=e^{-\ll}$.
%%%%%%%%%bbbbbbbbbbbbbbbbbbbbbbbbbbbbbbbbbbbbbb
\item[{ $b.$}]Our approach is similar to  that in the previous part. Here,  %$-\log n +\lnln \le c\le -\log\log n$ we have 
$e^{-c}\le 1$. Thus,
\begin{equation}\label{300.0}
    \begin{aligned}  \sum_{i=0}^{n-1}\left({\ec\log n}/{n}\right)^{    \el\ql^i}
\le e^{-c(1-o(1))}\sum_{i=0}^{n-1}\left({\log n}/{n}\right)^{    \el\ql^i}.
    \end{aligned}
\end{equation} 
Consider the sum on the \RHS of (\ref{300.0}). Taking $$g(n)=n-\lemmafunctionG,$$
we have
\begin{equation}\label{300.1}
    \begin{aligned}  \sum_{i=0}^{n-1}\left(\frac{\log n}{n}\right)^{    \el\ql^i}=&
    \sum_{0\le i\le g(n)}\left(\frac{\log n}{n}\right)^{    \el\ql^i}
    +\sum_{g(n)<i< n-\ll}\left(\frac{\log n}{n}\right)^{    \el\ql^i}\\
    &+\sum_{n-\ll\le i\le n-1}\left(\frac{\log n}{n}\right)^{    \el\ql^i}\\
    \le& \:n\cdot\left({\log n}/{n}\right)^{    \el\ql^{g(n)}}\\
    &+\left(n-g(n)\right)\cdot\left({\log n}/{n}\right)^{\el\ql^{n-\ll}}+\ll.
    \end{aligned}
\end{equation}
We start with the logarithm of the exponent in the first addend on the \RHS of (\ref{300.1}), 
\begin{equation*}%\label{}
    \begin{aligned}   
    \log \left(    \el\ql^{g(n)}\right)
    &=\ll+g(n)\log\ql\\
    &=\ll+\left(n-\lemmafunctionG\right)\left(-\ll/n+O(1/n^2)\right)\\
    &=\lnln/\log n+O(1/n).
    \end{aligned}
\end{equation*}
Thus, 
\begin{equation*}%\label{300.2}
    \begin{aligned}   
    \el\ql^{g(n)}
    &=\exp\left(\lnln/\log n+O(1/n)\right)\\
    &=1+\lnln/\log n+O(\lln2/\log^2 n),
    \end{aligned}
\end{equation*}
and hence  
\begin{equation*}%\label{}
    \begin{aligned}   
    &\log \left(    n\cdot\left(\log n/n\right)^{    \el\ql^{g(n)}}\right)\\
    &\hspace{10mm}=\log n+\el\ql^{g(n)}\left(-\log n+\log\log n\right)\\
    &\hspace{10mm}=\log n+\left(1+\lnln/\log n+O(\lln2/\log^2 n)\right)\left(-\log n+\log\log n\right)\\
    &\hspace{10mm}=o(1).
    \end{aligned}
\end{equation*}
Thus,
\begin{equation}\label{300.3}
    \begin{aligned}   
    n\cdot\left(\log n/n\right)^{    \el\ql^{g(n)}}
    &=\exp\left(o(1)\right)
    =1+o(1).
    \end{aligned}
\end{equation}
For the logarithm of the exponent in the second addend on the \RHS of (\ref{300.1}):
\begin{equation*}%\label{}
    \begin{aligned}   
    \log \left(    \el\ql^{n-\ll}\right)
    &=\ll+(n-\ll)\left(-\ll/n-\tfrac12\ll^2/n^2+O(1/n^3)\right)\\
    &=\tfrac12\ll^2/ n+O(1/n^2)>0.
    \end{aligned}
\end{equation*}
This implies that 
\begin{equation*}%\label{}
    \begin{aligned}   
        \el\ql^{n-\ll}> 1.
    \end{aligned}
\end{equation*}
Thus,
\begin{equation}\label{300.4}
    \begin{aligned}   \left(n-g(n)\right)\cdot\left({\log n}/{n}\right)^{\el\ql^{n-\ll}}
    &\le\left(n-g(n)\right)\cdot{\log n}/{n}\\
    &=\lemmafunctionG\cdot\frac{\log n}{n}
    =\lnln/\ll.
    \end{aligned}
\end{equation}
Altogether by (\ref{300.1})-(\ref{300.4})
\begin{equation}\label{300.5}
    \begin{aligned}  &\sum_{i=0}^{n-1}\left({\log n}/{n}\right)^{    \el\ql^i}\\
    %&\hspace{10mm}\ge\left(1-\frac{1}{\log n}\right)\cdot\frac{\log^{e^{\ll}} n}{n^{e^{\ll}-1}}+e^{-\ll}\left(1+\ll\lnln/\log n+O(1/\log n)\right)\\
    %&\hspace{10mm}\ge e^{-\ll}\\
    &\hspace{10mm}\le1+o(1)+\lnln/\ll+\ll\le \tfrac{2}{\ll}\lnln.
    \end{aligned}
\end{equation}
The claim follows from (\ref{300.0}) and (\ref{300.5}) with $\thB= 2/\ll$.

\end{description}

\qedForProof
\end{Proof}

\bigskip

%2222222222222222222222222222222222222222222222
\begin{Proof}{\ of Theorem \ref{Expectation}}

\begin{description}
\item[{ $a.$}]
Consider the  coupling of $D$ and $T$ presented in the proof of Lemma \ref{lemmaF_DF_T}.
Let~$\widetilde{T}_i$ be the time between the arrival of the $(i-1)$-st coupon (real or dummy) and that of the $i$-th coupon (real or dummy) in that process, $i\ge 1$ (where we agree that ``$0$-th coupon" arrives at time~$0$). 
Thus, $(\widetilde{T}_i)_{i=1}^{\infty}$ is a sequence of independent $\textnormal{Exp}(1)$-distributed variables.

We may  write:
\begin{equation*}
    \begin{aligned}
        \widetilde{T}=\widetilde{T}_1+\cdots +\widetilde{T}_{\widetilde{D}}.
    \end{aligned}
\end{equation*}
Therefore, as $\widetilde{T}_i$ are $\Exp(1)$-distributed, $1\le i\le \widetilde{D}$, 
\begin{align*}
    E\left(\widetilde{T}\vert\widetilde{D}\right)
=E\left(\widetilde{T}_1+\cdots +\widetilde{T}_{\widetilde{D}}\vert\widetilde{D}\right)
=\sum_{i=1}^{\widetilde{D}}E\left(\widetilde{T}_i\vert\widetilde{D}\right)
=\widetilde{D}.
\end{align*}
Consequently:
\begin{equation*}
    \begin{aligned}
        E({T})=E(\widetilde{T})
        =E\left(E\left(\widetilde{T}\vert \widetilde{D}\right)\right)
        =E\left( \widetilde{D}\right)=E(D).
    \end{aligned}
\end{equation*}

\item[{ $b.$}] 
Note first that  (\ref{eq-expectation}) cannot possibly follow from Theorem~\ref{ConvInDist'} %**
by itself. Thus, we return to the proof of that theorem and continue from there.

As $T$ is positive,
\begin{equation*}
    E(T)=\int\limits_{0}^{\infty}\left(1-F_T(t)\right)dt.
\end{equation*}
Changing variables,
$t=\tfrac{ne^{\lambda}}{\lambda} (\log n-\log \log n+c)$, 
% $$c=\frac{t\lambda}{e^{\lambda}n} -\log n+\log \log n,$$ or, equivalently, $$t=\frac{ne^{\lambda}}{\lambda} (\log n-\log \log n+c), \qquad dt=\frac{ne^{\lambda}}{\lambda} dc.$$ 
we obtain:
\begin{equation}\label{70.1}
    E(T)=\frac{ne^{\lambda}}{\lambda}\int\limits_{-\log n+\log \log n}^{\infty}\left(1-F_T\left(\tfrac{ne^{\lambda}}{\lambda} (\log n-\log \log n+c)\right)\right)dc.
\end{equation}
Denote $\ell_n=\log n-\log \log n$, $n>1$.
We will estimate  the integral on the \RHS of (\ref{70.1}) by splitting the interval $[-\ell_n,\infty)$ into three sub-intervals:~$\left[-\ell_n,-\aA\lnlnln\right]$, $\left[-\aA\lnlnln,\aB\lnln\right]$, and $\left[\aB\lnln,\infty\right]$. Denote by $I_j$ the integral on the $j$-th sub-interval, $1\le j\le 3$. We  estimate each $I_j$ separately.

We start with $I_1$.  By \Eq{1} and \Eq{2},  % WE GET $x>0$ $1-x\le e^{-x}$. I SUPPOSE THERE IS A STRAIGHT FORWARD REFERENCE FOR THIS, THOUGH IT IS A RESULT OF THE LEMMA ), 
using Lemma \ref{LemmaForExpectation}.a for large $n$ and some $\thA>0$,
\begin{equation*}%\label{7.thm2}
\begin{aligned}
F_T\left(t\right)
&=\prod_{i=0}^{n-1}\left(1-\left({\ec\log n}/{n}\right)^{
    \el\ql^i}\right)
    \le&\exp\left({-
\sum_{i=0}^{n-1}
\left({\ec\log n}/{n}\right)^{    \el\ql^i}}\right)\\
&\le \exp\left(-e^{-c(1-o(1))}\cdot\thA\right).
\end{aligned}
\end{equation*}
Thus, 
\begin{equation}\label{I1.1}
    \begin{aligned}
       I_1&=\int\limits_{-\ell_n}^{-\aA\lnlnln}\left(1-F_T\left(\tfrac{ne^{\lambda}}{\lambda} (\log n-\log \log n+c)\right)\right)dc\\
       &\ge\int\limits_{-\ell_n}^{-\aA\lnlnln}\left(1-\exp\left(-e^{-c(1-o(1))}\cdot\thA\right)\right)dc\\
       &=\log n-\lnln-\aA\lnlnln-\int\limits_{-\ell_n}^{-\aA\lnlnln}\exp\left(-e^{-c(1-o(1))}\cdot\thA\right)dc.
    \end{aligned}
\end{equation}
Consider the integral on the \RHS of (\ref{I1.1}).
Note that, as $ c\le -\aA\lnlnln<0$, we have $e^{-2c/3}>1$. %$\log^{\aA}\log n\le e^{-c(1-o(1))}\le \left(n/\log n\right)^{1-o(1)}$.
Thus, for large $n$
%As $e^{-c(1-o(1))}\ge 1$,
\begin{equation}\label{I1.2}
    \begin{aligned}
       0&\le\int\limits_{-\ell_n}^{-\aA\lnlnln}\exp\left(-e^{-c(1-o(1))}\cdot\thA\right)dc\\
       &\le \frac{1}{\thA}\int\limits_{-\ell_n}^{-\aA\lnlnln}\thA e^{{-2c/3}}\exp\left(-e^{-2c/3}\cdot\thA\right)dc\\%note that \exp(-e^{-c(1-o(1))}\thA)\le \exp(-e^{-c/2}\thA)
       &= \left[\frac{3\exp\left(-e^{-2c/3}\cdot\thA\right)}{2\thA}\right]
       _{-\ell_n}^{-\aA\lnlnln}\\
       &= O\left(\exp\left(-\thA(\log\log n)^{4/3}\right)\right)=O\left({1}/{\log n}\right).
    \end{aligned}
\end{equation}
As $I_1\le \log n-\lnln-\aA\lnlnln$,  by (\ref{I1.1}) and (\ref{I1.2})
\begin{equation}\label{I1}
    \begin{aligned}
       I_1=\log n-\lnln-\aA\lnlnln+O\left({1}/{\log n}\right).
    \end{aligned}
\end{equation}

Skip to $I_3$.
For $c\ge\aB \log\log n$, by (\ref{7}) and Lemma~\ref{LemmaForExpectation}.b for  some $\thB>0$,
\begin{equation*}%\label{7.largec.1}
\begin{aligned}
   &F_T\left(\tfrac{ne^{\lambda}}{\lambda} (\log n-\log \log n+c)\right)\\
   &\hspace{10mm}=\exp\left({-
\sum_{i=0}^{n-1}
\left({\ec\log n}/{n}\right)^{    \el\ql^i}}\right)
   +O\left(\frac{e^{-2c+o(1)}\log^2 n}{{n}^{
    1-o(1)}}\right)\\
   &\hspace{10mm}\ge \exp\left(-e^{-c(1-o(1))}\cdot\thB\lnln\right)
   +O\left({e^{-c}\log^2 n}/{n}^{
    1-o(1)}\right).
\end{aligned}
\end{equation*}
Thus, by a routine calculation:
\begin{equation}\label{I3}
    \begin{aligned}
       0\le I_3&=\int\limits_{\aB\log \log n}^{\infty}\left(1-F_T\left(\tfrac{ne^{\lambda}}{\lambda} (\log n-\log \log n+c)\right)\right)dc\\
       &\le\int\limits_{2\log \log n}^{\infty}\left(1-\exp\left(-e^{-c(1-o(1))}\cdot\thB\lnln\right)+O\left(\frac{e^{-c}\log^2 n}{{n}^{
    1-o(1)}}\right)\right)dc\\
       &=\int\limits_{2\log \log n}^{\infty}\left(1-\left(1+O\left(e^{-c(1-o(1))}\cdot\lnln\right)\right)+O\left(\frac{e^{-c}\log^2 n}{{n}^{
    1-o(1)}}\right)\right)dc\\
       &=\int\limits_{2\log \log n}^{\infty}O\left(e^{-c(1-o(1))}\cdot\lnln\right)dc=O\left({1}/{\log n}\right).
    \end{aligned}
\end{equation}
We will bound  $I_2$  from both sides, using the bounds obtained above on $F_T$. We start with an upper bound. By  (\ref{69.3}):
\begin{equation}\label{70.2}
\begin{aligned}
&F_T\left(\tfrac{ne^{\lambda}}{\lambda} (\log n-\log \log n+c)\right)\\
%&\hspace{1cm}\ge e^{-e^{-c}/\ll}\cdot\exp\left(-\frac{ e^{-c}}{\ll} \left(\lglgFlgly{} +\errLognly{}\right)\right)\\
&\hspace{1cm}\ge e^{-e^{-c}/\ll}\cdot\left( 1-\frac{ e^{-c}}{\ll} \left(\lglgFlgly{}-{ce^{-\la}}/{\log n} +\errLognly{}\right)\right).
\end{aligned}
\end{equation}
 By \Eq{70.1} and \Eq{70.2},
\begin{equation}\label{70.3}
\begin{aligned}
I_2
\le& \int\limits_{-2\lnlnln}^{2\lnln}\left(1-e^{-e^{-c}/\lambda} \left(1-\frac{ e^{-c}}{\ll} \left(\lglgFlgst{}-\frac{ce^{-\la}}{\log n} +O\left(\frac{1}{\log n}\right)\right)\right)\right)dc\\
=&\int\limits_{-2\lnlnln}^{2\lnln}\left(1-e^{-\ec/\ll}\right)dc
+\left(\lglgFlgst{}+O\left(\frac{1}{\log n}\right)\right)\int\limits_{-2\log\log\log n}^{2\lnln}\frac{e^{-c}}{\ll} e^{-\ec/\ll}dc\\
&-\frac{e^{-\la}}{\log n}\int\limits_{-2\lnlnln}^{2\lnln}c\cdot\frac{e^{-c}}{\ll} e^{-\ec/\ll}dc.  
\end{aligned}
\end{equation}
We may rewrite \Eq{70.3} in the form
\begin{equation}\label{72.1a}
\begin{aligned}
I_2\le& \int\limits_{-2\lnlnln}^{2\lnln}\left(1-F_G(c)\right)dc
+\left(\lglgFlgst{}+O\left(\frac{1}{\log n}\right)\right)\int\limits_{-2\lnlnln}^{2\lnln}f_G(c)dc\\
&-\frac{e^{-\la}}{\log n}\int\limits_{-2\lnlnln}^{2\lnln}c\cdot \frac{e^{-c}}{\ll} e^{-\ec/\ll}dc\\
%=& \int\limits_{-2\lnlnln}^{2\lnln}\left(1-F_G(c)\right)dc-\frac{e^{-\la}}{\log n}\int\limits_{-2\lnlnln}^{2\lnln}c\cdot \frac{e^{-c}}{\ll} e^{-\ec/\ll}dc+O\left(\frac{\lnln}{\log n}\right)\\
=&\int\limits_{-2\lnlnln}^{2\lnln}\left(1-F_G(c)\right)dc
+O\left({\lnln}/{\log n}\right),
\end{aligned}
\end{equation}
where  $G\sim \textnormal{Gumbel}(-\log \lambda, 1)$, and  $F_G$ and $f_G$ are its distribution and density functions, respectively. %\Eq{gumbel}, 
%\begin{equation}\label{G}F_G(c)=e^{-e^{-c}/\lambda},\qquad f_G(c)=\tfrac{1}{\lambda}e^{-c}\cdot e^{-e^{-c}/\lambda},\qquad -\infty<c<\infty.\end{equation}
%explanation on the deletion of the second row into O(\log\log n/\log n):Integrating by parts  the integral in the second addend on the \RHS of (\ref{72.1a}) we estimate it to be  of size magnitude  $\lnln$. 
%\begin{equation*}    \begin{aligned}       \int\limits_{-2\lnlnln}^{2\lnln}c\cdot f_G(c)dc=&\Big[c\cdot F_G(c)\Big]  _{-2\lnlnln}^{2\lnln}  -\int\limits_{-2\lnlnln}^{2\lnln}F_G(c)dc=O(\lnln).  \end{aligned}\end{equation*}
%LLLLLOWER------

For the lower bound on $I_2$, we proceed similarly. By (\ref{1a}) and (\ref{70.0}):
\begin{equation}\label{70.2b}
\begin{aligned}
&F_T\left(\tfrac{ne^{\lambda}}{\lambda} (\log n-\log \log n+c)\right)\\
&\hspace{1cm}\le e^{-e^{-c}/\ll}\cdot\left(1-\frac{ e^{-c}}{\ll} \left((1-\la e^{-\la})\lglgFlgst-\frac{ce^{-\la}}{\log n}+O\left(\frac{1}{\log n}\right)\right)\right).
\end{aligned}
\end{equation}
By \Eq{70.1}
 and \Eq{70.2b},
\begin{equation}\label{72.1b}
\begin{aligned}
I_2
\ge& \int\limits_{-2\lnlnln}^{2\lnln}\left(1-e^{\tfrac{e^{-c}}{\ll}} \left(1-\tfrac{ e^{-c}}{\ll} \left((1-\la e^{-\la})\lglgFlgtst-\tfrac{ce^{-\la}}{\log n}+\errLognst\right)\right)\right)dc\\
=&\int\limits_{-2\lnlnln}^{2\lnln}\left(1-F_G(c)\right)dc
+\left((1-\la e^{-\la})\lglgFlgtst+\errLognst\right)\int\limits_{-2\lnlnln}^{2\lnln}f_G(c)dc\\
&-\frac{e^{-\la}}{\log n}\int\limits_{-2\lnlnln}^{2\lnln}c\cdot \frac{e^{-c}}{\ll} e^{-\ec/\ll}dc\\
=&\int\limits_{-2\lnlnln}^{2\lnln}\left(1-F_G(c)\right)dc
+O\left(\lnln/\log n\right),
\end{aligned}
\end{equation}
Altogether, by (\ref{72.1a}) and (\ref{72.1b}), 
\begin{equation}\label{72.1.0}
\begin{aligned}
I_2=& \int\limits_{-2\lnlnln}^{2\lnln}\left(1-F_G(c)\right)dc
+O\left({\lnln}/{\log n}\right).
\end{aligned}
\end{equation}
We may calculate $I_2$ as follows:
\begin{equation}\label{72.1.0.1}
\begin{aligned}
I_2=& \int\limits_{-\ell_n}^{\infty}\left(1-F_G(c)\right)dc
-\int\limits_{-\ell_n}^{-2\lnlnln}\left(1-F_G(c)\right)dc\\
&-\int\limits_{2\lnln}^{\infty}\left(1-F_G(c)\right)dc
+O\left({\lnln}/{\log n}\right).
\end{aligned}
\end{equation}
We start with the second integral on the \RHS of (\ref{72.1.0.1}).
Going over the calculations in (\ref{I1.1})-(\ref{I1}), we notice  that, if we replaced $\exp\left(-e^{-c(1-o(1))}\thA\right)$ in (\ref{I1.1}) by $F_G(c)$, we would still get the same result as in (\ref{I1}). Namely,
\begin{equation}\label{72.1.0.I1}
\begin{aligned}
 I_1=\int\limits_{-\ell_n}^{-2\lnlnln}\left(1-F_G(c)\right)dc
 %-\log n+\lnln+O\left(\lnln/\log n\right)=
 +O\left({1}/{\log n}\right).
\end{aligned}
\end{equation}
Similarly, replacing $\exp\left(-e^{-c(1-o(1))}\thB\lnln\right)$ in  (\ref{I3}) by $F_G(c)$,  the third integral on the \RHS of (\ref{72.1.0.1}) becomes $O(1/\log^2 n)$.
Thus,
\begin{equation}\label{72.1.0.I3}
\begin {aligned}
 I_3=\int\limits_{2\lnln}^{\infty}\left(1-F_G(c)\right)dc
 +O\left(1/\log n\right).
\end{aligned}
\end{equation}
By (\ref{70.1}) and (\ref{72.1.0})-(\ref{72.1.0.I3}) we have:
\begin{equation}\label{72.1}
\begin{aligned}
E(T)=\frac{n\el}{\ll}\left(I_1+I_2+I_3\right)
=\frac{n\el}{\ll} \int\limits_{-\ell_n}^{\infty}\left(1-F_G(c)\right)dc
+O\left({\lnln}/{\log n}\right).
\end{aligned}
\end{equation}
Consider the  integral on the \RHS of (\ref{72.1}).
For a variable $Y$ with distribution function $F_Y$ and density $f_Y$, by \citep[p.150, (6.3)]{Feller1970},
\begin{equation}\label{72.2}
\int\limits_0^{\infty}yf_Y(y)dy=\int\limits_0^{\infty}(1-F_Y(y))dy,
\end{equation}
and by \citep[p.150, (6.1)]{Feller1970}, for  $b>0$,
\begin{equation}\label{73.1}
\begin{aligned}
    \int\limits_{-b}^0 yf_Y(y)dy
    &=-(-b)F_Y(-b)-\int\limits_{-b}^0 F_Y(y)dy\\
    &=-b\left(1-F_Y(-b)\right)+\int\limits_{-b}^0 \left(1-F_Y(y)\right)dy.
\end{aligned}
\end{equation}
Thus, by \Eq{72.2} and \Eq{73.1},
\begin{equation}\label{73.1a}
\int\limits_{-b}^{\infty}\left(1-F_G(c)\right)dc=\int\limits_{-b}^{\infty}c\cdot f_G(c)dc+b\left(1-F_G(-b)\right).
\end{equation}
By \Eq{72.1} and \Eq{73.1a}, for $b=\ell_n$, 
\begin{equation}\label{73.2}
\begin{aligned}
E(T)&= \frac{ne^{\lambda}}{\lambda}\int\limits_{-\ell_n}^{\infty}c\cdot f_G(c)dc+\frac{ne^{\lambda}}{\lambda}\cdot\ell_n\left(1-F_G(-\ell_n)\right)+O\left(\frac{n\lnln}{\log n}\right).    
\end{aligned}
\end{equation}
Now:
\begin{equation}\label{74.1}
\begin{aligned}
F_G(-\ell_n)
&=\exp\left(-e^{-(-\ell_n)}/\ll\right)
=\exp\left(-e^{\log n-\lnln}/\ll\right)
=\exp\left(-{n}/({\ll\log n})\right).
\end{aligned}
\end{equation}
Using integration by parts,
\begin{equation}\label{74.2}
\begin{aligned}
\int\limits_{-\ell_n}^{\infty}c\cdot f_G(c)dc
=\:&E(G)-\int\limits_{-\infty}^{-\ell_n}c\cdot f_G(c)dc
=E(G)-\Big[c\cdot F_G(c)\Big]_{-\infty}^{-\ell_n}+\int\limits_{-\infty}^{-\ell_n}F_G(c)dc. 
\end{aligned}
\end{equation}
Recall that, if $X\sim\textnormal{Gumbel}(\mu,\beta)$,  then $E(X)=\mu+\beta\gamma$.  Thus,
\begin{equation}\label{74.3}
E(G)=-\log\ll+\gamma.
\end{equation}
Also,
\begin{equation}\label{74.35}
\lim_{c\to -\infty} c\cdot F_G(c)
=\lim_{c\to -\infty} c\exp\left(-e^{-c}/\ll\right)=0.
\end{equation}
Thus, by (\ref{74.1}) and (\ref{74.35}):
\begin{equation}\label{74.4}
\Big[c\cdot F_G(c)\Big]_{-\infty}^{-\ell_n}
=-\ell_n F_G(-\ell_n)
=\left(-\log n+\lnln\right) \exp\left(-{n}/({\ll\log n})\right).
\end{equation}
Consider the third addend on the \RHS of (\ref{74.2}). By (\ref{74.1}),
\begin{equation}\label{75.1}
\begin{aligned}
0&\le\int\limits_{-\infty}^{-\ell_n}F_G(c)dc
=\int\limits_{-\infty}^{-\ell_n}\exp\left(e^{-c}/\ll\right)dc
\le\ll\int\limits_{-\infty}^{-\ell_n}\frac{e^{-c}}{\ll}\cdot\exp\left(e^{-c}/\ll\right)dc\\
&=\ll\int\limits_{-\infty}^{-\ell_n}f_G(c)dc=\ll F_G(-\ell_n)=\ll\exp\left(-{n}/({\ll\log n})\right). 
\end{aligned}
\end{equation}
Altogether, by (\ref{74.2})-(\ref{75.1}):
\begin{equation}\label{75.2}
\begin{aligned}
&\int\limits_{-\ell_n}^{\infty}c\cdot f_G(c)dc\\
&\hspace{1cm}=-\log\ll+\gamma-\left(-\log n+\lnln\right) \exp\left(-\tfrac{n}{\ll\log n}\right)
%&\hspace{225pt}
+O\left(\exp\left(-\tfrac{n}{\ll\log n}\right)\right)\\
&\hspace{1cm}=-\log\ll+\gamma
+O\left(\log n\cdot\exp\left(-{n}/({\ll\log n})\right)\right).
\end{aligned}
\end{equation}
By (\ref{73.2}), (\ref{74.1}) and (\ref{75.2}),
\begin{equation}\label{75.3}
\begin{aligned}
E(T)=
\:&\frac{n\el}{\ll}
\left(-\log\ll+\gamma
+O\left(\log n\cdot\exp\left(-{n}/({\ll\log n})\right)\right)\right) \\
&+\frac{n\el}{\ll}\cdot\ell_n\left(1-\exp\left(-{n}/({\ll\log n})\right)\right)+O\left({n\lnln}/{\log n}\right)\\
=&\:\frac{n\el}{\ll}\left(\log n-\lnln-\log\ll+\gamma
+O\left(\lnln/\log n\right)\right).    
\end{aligned}
\end{equation}
\end{description}

\qedForProof
\end{Proof}

\bibliography{refs}
\bibliographystyle{abbrv}

\end{document}